\documentclass[12pt, makeindex]{amsbook}
\usepackage[english]{babel}
\usepackage{amsmath,amsfonts}
\usepackage{amssymb}
\usepackage{amscd}
\usepackage{graphics}
\usepackage{graphicx}
\usepackage[all]{xy}
\usepackage{pstcol}

\newtheorem{theorem}{Theorem}[section]
\newtheorem{proposition}[theorem]{Proposition}
\newtheorem{lemma}[theorem]{Lemma}
\newtheorem{corollary}[theorem]{Corollairy}
\newtheorem{definition}[theorem]{Definition}
\newtheorem{example}[theorem]{Example}
\newtheorem{examples}[theorem]{Examples}
\newtheorem{remark}[theorem]{Remark}

\def\Tr{\mathop{\mathbf{Tr}}}
\def\rank{\mathop{\mathrm{rank}}}
\def\tq{\hskip0.5mm \mid \hskip0.5mm} 

\def\c{\hbox{\tt c}}
\def\C{\hbox{\tt C}}

\def\tr{\hbox{\rm tr}}

\def\ie{{\it i.e}} 

\let\ord=\preccurlyeq
\def\tr{\mathop{\mathbf{tr}}}
\def\det{\mathop{\mathbf{det}}}
\def\Det{\mathop{\mathrm{Det}}}

\makeindex
\begin{document}

\title{\Huge \vskip6cm Jordan algebras, geometry of Hermitian symmetric spaces and non-commutative Hardy spaces \vskip2cm }
\author{Khalid Koufany }
\author{}

\address{
Institut {\'E}lie Cartan, 
Universit{\'e} Henri Poincar{\'e}, (Nancy 1)
B.P. 239, F-54506 Vand{\oe}uvre-l{\`e}s-Nancy cedex, France}
\email{Khalid.Koufany@iecn.u-nancy.fr}

\maketitle


These notes were written following lectures I had the pleasure of giving on this subject at Keio University, during November and December 2004.\\

The first part is  about  new applications of Jordan algebras to the geometry of Hermitian symmetric spaces and to causal semi-simple symmetric spaces of Cayley type. 

The second part will present new contributions for studying  (non commutative) Hardy spaces of holomorphic functions on Lie semi-groups 
which is a  part of the so called {\it Gelfand-Gindikin program}.\\

I wish to acknowledge the support of  ``{\it The 21st Century CEO Program
  at Keio}''. \\
I wish to express my gratitude to Prof. T. Kawazoe and Prof. Y. Maeda for inviting me to give these lectures. \\
I also thank Prof. T. Oshima of Tokyo Metropolitan University, Prof. T. Nomura of Kyoto University and  Prof. H. Ishi of Yokohama city University, where a part of the lectures were also given. \\
Finally I am grateful to Prof. T. Kawazoe for inviting me to Keio University, for his kindness and hospitality.

\tableofcontents




Let $G/K$ be a Riemannian symmetric space of non compact type. Suppose that $G/K$ is an
irreducible  Hermitian symmetric space, {\it i.e.\/} the center
$\frak{z}(\frak k)$ of  $\frak k$ (the Lie
algebra  of $K$) is non trivial.\\

 There are two essentially equivalent ways to treat Hermitian symmetric spaces.

 The first one uses the theory of semi-simple Lie groups; it was in this way that the basic facts of the theory were established by {\'E}lie Cartan in the 1930's and by Harish-Chandra in the 1950's.

 The second way avoids semi-simple Lie groups theory and uses Jordan
 algebra and triple systems; it is essentially due to Koecher and his
 school in the 1980's. \\

In our presentation we will mainly use the Jordan theory to point out some
development in the geometry and analysis on Hermitian symmetric spaces.\\

Let $\theta$ be the Cartan involution of the Hermitian symmetric space
$G/K$. There exists a non trivial involution $\tau$ of the group $G$ which
commutes with the given Cartan involution. Let $H$ be the $\tau-$fixed
points of $G$. Then the symmetric space $G/H$ belongs to an important class
of non Riemannian symmetric spaces, namely, the class of causal symmetric
spaces introduced by {\'O}lafsson and {\O}rsted in the 1990's.\\

 Suppose  $G/K$ of tube type, holomorphically equivalent, via the Cayley
 transform, to the tube $V+i\Omega$, where $V$ is a Euclidean Jordan
 algebra and $\Omega$ its symmetric cone. We prove (see section
 \ref{section4}) that in this case $G/H$ is a causal symmetric space of
 Cayley type. In section \ref{section5} we give a conformal compactification of such spaces. We
 also investigate the semigroup associated with the order on $G/H$, see
 section \ref{section6}. We
 prove that it is related to the semigroup $S_\Omega$ of compressions of
 $\Omega$. Each element of $S_\Omega$ is a contraction for the Riemannian
 metric as well as for the Hilbert metric of
 $\Omega$.\\

Let $S$ be the Shilov boundary of $G/K$. In section \ref{section6}, we study the causal structure
of $S$. This boundary has many geometric invariants, in particular the
transversality index (see
section \ref{section7}) and
the triple Malsov index (see section \ref{section8}). The universal covering of $S$ is needed to
study other  geometric invariants introduced by Souriau and Arnold-Leray
in the Lagrangian case. In section \ref{section9} we give an explicit construction of this universal
covering. We also generalize the Souriau index (see section
\ref{section10}) and the Arnold-Leray index (see section
\ref{section11}). Finally, in section \ref{extra_section}, we use the
Souriau index to 
generalize  of the Poincar{\'e} rotation number.\\

A further interesting development is the so-called Gelfand-Gindikin program. In 1977
 Gelfand and Gindikin proposed a new approach to study the Plancherel
 formula for semi-simple Lie groups. The idea is to consider functions on
 $G$ as boundary values of holomorphic functions of a domain (Lie
 semi-groups) in the
 complexification of $G$ and to study the action of $G$ on these
 holomorphic functions. When $G$ is the group of holomorphic diffeomorphisms
 of a Hermitian symmetric space of tube
 type, the Gelfand-Gindikin program has been developed by Olshanski\u{\i} in several papers. This study is
related to harmonic analysis of bounded symmetric domains and the
decompositions of the Hardy spaces of Lie semi-groups, which involves the holomorphic discrete
series of representations of $G$.\\

We present in section \ref{section12} the theory of Olshanski\u{\i}. In
sections \ref{section15}, \ref{section16} and \ref{section17} we
develop this program for the groups $Sp(r,\mathbb{R})$, $SO^*(2\ell)$ and
$U(p,q)$. We give a new approach of studying Hardy spaces on Lie
semi-groups. More precisely, we introduce a 
new Cayley transform, which allows us to compare the classical Hardy space and
the Hardy space on the Lie semi-group.

\newpage

\part{Cayley type  symmetric spaces, transversality and the Maslov index}

\section{Causal symmetric spaces}\label{section1}

{\it In this section we will recall the notion of causal symmetric
  spaces introduced by {\'O}lafsson and {\O}rsted and give some examples. For
noncompactly causal symmetric spaces, we will introduce the corresponding
Olshanski\u{\i} semigroup.}

\subsection{Causal structures}
Let $V$ be a vector space. A subset $C$ of $V$ is called a  {\it causal cone} \index{Cone!causal}  if it is non-zero closed convex and proper (i.e. $C\cap -C=\{0\}$) cone.\\
Let $\mathcal{M}$ be a $n-$dimensional manifold.  A {\it causal structure}\index{Structure!causal} on $\mathcal{M}$ is an assignment, $x\mapsto C_x$, to each point $x$ of $\mathcal{M}$ a causal cone $C_x$ in the tangent space $T_x(\mathcal{M})$ of $\mathcal{M}$ at $x$  such that  $C_x$ depends smoothly on  $x$.\\
A  $\mathcal{C}^1$ curve $\gamma :[\alpha,\beta]\to \mathcal{M}$ is said to {\it causal curve} (resp; {\it anti-causal curve})\index{Causal curve} if $\dot{\gamma}(t)\in C_{\gamma(t)}$  (resp. $\dot{\gamma}(t)\in -C_{\gamma(t)}$) for all $t$. If there is no non-trivial closed causal curves, the causal structure of $\mathcal{M}$ is said to {\it global}. We can then define a partial order $\preccurlyeq$ on $\mathcal{M}$ : 
 $$x\preccurlyeq y \; \textrm{if there exists a causal curve} \; \gamma : [\alpha,\beta]\to \mathcal{M}, \; \gamma(\alpha)=x,\; \gamma(\beta)=y. $$
If $\mathcal{M}=G/H$ is a homogeneous space, where $G$ is a Lie group and $H$ is a closed subgroup of $G$, then  the causal structure of $\mathcal{M}$ is said to be  $G-${\it invariant} if, for every $g\in G$,
$$C_{g\cdot x}=Dg(x)(C_x),$$
where $Dg(x)$ is the derivative of $g$ at $x$.
Let $x_o=eH$ be the base point. An invariant  causal structure on $\mathcal{M}=G/H$ is determined by a  causal cone $C_o\subset T_{x_o}(\mathcal{M})$ which is invariant under the action of $H$.
\subsection{Causal symmetric spaces}
Suppose that  $\mathcal{M}=G/H$ is a symmetric space  :  there exists an involution $\sigma$ of $G$ such that $(G^\sigma)^\circ\subset H\subset G^\sigma$, where $(G^\sigma)^\circ$ is the connected component of $G^\sigma$, the subgroup  of  fixed points of $\sigma$ in $G$. \\
 Let $\frak g$ be the Lie algebra of $G$. Put $\frak h=\frak g(+1,\sigma)$ and $\frak q=\frak g(-1,\sigma)$ the eigenspaces of $\sigma$. Then $\frak h$ is the Lie algebra of $H$ and $\frak g=\frak h\oplus \frak q$. The tangent space of $\mathcal{M}$ at the point $x_o$ can be identified with $\frak q$. With this identification, the derivative $Dh(x_o)$, $h\in H$ corresponds to $\textrm{Ad}(h)$. Therefore an invariant causal structure on $\mathcal{M}$ is determined by a  causal cone $C\subset \frak q$ which is $\textrm{Ad}(H)-$invariant.\\
 
 Suppose  $G$ semi-simple and has a finite center and suppose that the pair $(\frak g,\frak h)$ is irreducible ( i.e. there is no non-trivial ideal in $\frak g$, invariant by $\sigma$). In this case there exists a Cartan involution $\theta$ of $G$ such that $\sigma\theta=\theta\sigma$. Let $K=G^\theta$, then $K$ is a maximal compact subgroup of $G$.  Let $\frak k=\frak g(+1,\theta)$ and $\frak p=\frak g(-1,\theta)$.  Then  $\frak k$ is the Lie algebra of $K$ and   $\frak g=\frak k\oplus \frak p$ is the Cartan decomposition of $\frak g$.  Moreover we have
 $$ 
 \frak g=(\frak h\cap{\frak k})\oplus (\frak q\cap{\frak k})\oplus(\frak h\cap{\frak p)}\oplus(\frak q\cap{\frak p}).
 $$
 Let $\textrm{Cone}_H(\frak q)$ be the set of $\textrm{Ad}(H)-$invariant causal cones in $\frak q$.
 \begin{definition}
 Let  $\mathcal{M}=G/H$  be an irreducible non-Riemannian semi-simple symmetric space. Then we call $\mathcal{M}=G/H$ 
 \begin{enumerate}
\item[( \rm{CC})]   compactly causal symmetric space,\index{Symmetric space!compactly causal}  if there exists $C\in \textrm{Cone}_H(\frak q)$ such that $C^\circ\cap\frak k\not=\emptyset$.
\item[( \rm{NCC})]   a non-compactly causal symmetric space,\index{Symmetric space!non-compactly causal} if there exists $C\in \textrm{Cone}_H(\frak q)$ such that $C^\circ\cap\frak p\not=\emptyset$.
\item[( \rm{CT})]   symmetric space of Cayley type,\index{Symmetric space!Cayley type} if both (CC) and (NCC) hold.
\item[( \rm{CAU})]  causal symmetric space \index{Symmetric space!causal}if either (CC) or (NCC) holds.
\end{enumerate}
 \end{definition}
 Let 
 $$\frak q^{H\cap K}=\{X\in \frak q \; ,\; \forall k\in H\cap K \; :\; \textrm{Ad}(k)X=X\}.$$
 Then, there exists an $\textrm{Ad}(H)-$invariant causal cone in $\frak q$ if and only if $\frak q^{H\cap K}\not=\{0\}$. 
 \begin{proposition} Let $\mathcal{M}=G/H$ be an irreducible symmetric space.
 \begin{enumerate}
 \item $\mathcal{M}$ is compactly causal if and only if $\frak q^{H\cap K}\cap \frak k\not =\{0\}$.
 \item $\mathcal{M}$ is non-compactly causal if and only if $\frak q^{H\cap K}\cap \frak p\not =\{0\}$.
 \end{enumerate}
 \end{proposition}
\begin{examples} {\rm We give some examples of causal symmetric spaces.
  \begin{enumerate}
    \item[(CC)]
      \begin{enumerate}
      \item The group case. Let $G_1$ be a semi-simple Lie group. Let $G=G_1\times G_1$ and define $\sigma(a,b)=(b,a)$. Then $H=\{(a,a)\; ;\; a\in G_1\}\simeq G_1$. The symmetric space $\mathcal{M}=G/H\simeq G_1$ has an invariant causal structure if and only if the Lie algebra $\frak g_1$ of $G_1$ is of Hermitian type. In this case $\mathcal{M}$ is a compactly causal symmetric space.
      \item The hyperboloid $Q^n_+=\{v\in\mathbb{R}^{n+1} \; ;\; v_1^2+v_2^2-v_3^2-\ldots-v_{n+1}^2=1\}$.  We have  
$$Q^{1,n}_+\simeq SO_o(2,n)/SO_o(1,n)$$  
and it is a symmetric space where the involution  
  $\sigma$ is the conjugation by the matrix $$J_{1,n}=\begin{pmatrix}1&0\\0&-I_n\end{pmatrix}.$$ 
In this case $\frak q$ is the set of all matrices
$$q_+(v)=\begin{pmatrix}0&-v^t\\ v& 0\end{pmatrix}\; ;\; v\in\mathbb{R}^{n+1}$$
in which $SO_0(1,n)$ acts as usual. Let 
$$C_+=\{q_+(v)\in\frak q \; ; \; v_1^2-v_2^2-\ldots-1\geq 0, v_1\geq 0\}.$$
Then $C_+$ is a causal cone in $\frak q$ invariant under the  group $SO_0(1,n)$ and $\text{Cone}_{SO_0(1,n)}(\frak q)=\{C_+,-C_+\}$. In particular $SO_o(2,n)/SO_o(1,n)$ carries a compactly causal structure.
       \end{enumerate}
\item[(NCC)]
\begin{enumerate}
    \item The Ol'shanski{\u\i} symmetric spaces.\index{Symmetric space!Ol'shanski{\u\i}}
       Let $\frak g_1$ be a semi-simple Hermitian Lie algebra. Let $\frak g=(\frak g_1)_\mathbb{C}$ and let $\sigma$ be the complex  conjugation. let $G_\mathbb{C}$ be a complex analytic group with Lie algebra $\frak g$. Let $G_1$ be the analytic subgroup of $G_\mathbb{C}$ corresponding to $\frak g_1$. Then $G_1=(G^\sigma_\mathbb{C})^\circ$ and $\mathcal{M}=G_\mathbb{C}/G_1$ is a symmetric space. In this case $\frak q=i\frak g_1$ and $\mathcal{M}$ is non-compactly causal symmetric space.
       \item The hyperboloid $Q^n_-=\{v\in\mathbb{R}^{n+1} \; ;\; v_1^2-v_2^2-\ldots-v_{n+1}^2=-1\}$.  We have
       $$Q^n_-\simeq SO_0(1,n)/SO_0(1,n-1)$$
       and it is a symmetric space. In this case $\frak q$ is the set
      $$q_-(w)=\begin{pmatrix}0&w\\ w^t& 0\end{pmatrix}\; ;\; w\in\mathbb{R}^{n+1}$$ 
      in which the group $SO_0(1,n)$ acts. Let
      $$C_-=\{q_-(w)\in\frak q \; ; \; w_1^2-w_2^2-\ldots-1\geq 0, w_1\geq 0\}.$$
      Then $C_-$ is a causal cone in $\frak q$, invariant under the group $SO_0(1,n)$. Furthermore, $\text{Cone}_{SO_0(1,n)}(\frak q)=\{C_-,-C_-\}$ and $SO_0(1,n+1)/SO_0(1,n)$ is a non-compactly causal symmetric space.
      
    \end{enumerate}
   \item[(TC)]
     \begin{enumerate}
     \item[] 
      The hyperboloid of one sheet in $\mathbb{R}^3$, $\mathrm{SO}_o(2,1)/\mathrm{SO}_o(1,1)$ is a symmetric space of Cayley type. One can realize it as the off-diagonal subset of $S^1\times S^1$, where $S^1$ is the unit circle.
$$\mathcal{M}\simeq\{(u,v)\in S^1\times S^1 \; ;\; u\not=v\}.$$
We will see, in section \ref{sec_cayley_type},  that this happens for any symmetric space of Cayley type.

       \end{enumerate}
    \end{enumerate}
}
\end{examples}

If $\mathcal{M}=G/H$ is noncompactly causal symmetric space, then there exists a causal cone $C\in \textrm{Cone}_H(\frak q)$ such that $C^\circ \cap \frak p\not=\emptyset$ and $C\cap \frak k=\{0\}$. The causal structure is then global and we can define and order $\preccurlyeq$ on $\mathcal{M}$. The set
    $$S_\preccurlyeq =\{g\in G \; ;\; x_o\preccurlyeq g\cdot x_o\}$$
      is a closed semigroup  called the {\it Ol'shanski{\u \i} semigroup}.\index{Semi-group!Ol'shanski{\u \i}}  One has the  {\it the Olshanski{\u \i} decomposition}\index{Decomposition!Ol'shanski{\u \i}}
    $$S_\preccurlyeq =H\exp(C).$$


\newpage
 \section{Jordan algebras and symmetric cones}\label{section2}
{\it In this section we recall the notion of Euclidean Jordan algebras and
  fix notations. Our presentation is mainly based on \cite{Faraut-Koranyi}}.\\

Let $V$ be a {\it Euclidean Jordan algebra} \index{Jordan!algebra} with identity element $e$ and of dimension $n$. This means that $V$ is a $n-$dimensional Euclidean vector space equipped with a bilinear product such that
\begin{equation*}
xy=yx,
\end{equation*}
\begin{equation*}
x^2(xy)=x(x^2y),
\end{equation*} 
\begin{equation}\label{assoc}
(xy|z)=(x|yz).
\end{equation}
For $x\in V$, denote by $L(x)$ the linear operator defined by $y\mapsto L(x)y=xy$, and introduce the {\it quadratic representation} $P(\cdot)$ and the ''square'' operator $\Box$, defined by
$$P(x)=2L(x)^2-L(x^2),\;\;\; x\Box y=L(xy)+[L(x),L(y)],$$
where the brackets denote the commutator.
For simplicity we assume $V$ to be simple. In other words, there is non non-trivial ideal in $V$. 
An element $x$ is said to be {\it invertible} if there exists an element $y\in\mathbb{R}[x]$ such that $xy=e$. Since $\mathbb{R}[x]$ is associative, $y$ is unique. It is called the {\it inverse} of $x$ and is denoted by $y=x^{-1}$. Let $V^\times$ the set of invertible elements of  $V$. The connected component of the unit $e$ in $V^\times$ is the set  $\Omega$ of squares,
$$\Omega=\{x^2 \; ;\; x\in V^\times\}.$$
The set $\Omega$ is open, convex, proper, generating, symmetric, homogeneous cone. This is the {\it symmetric cone}\index{Cone!symmetric} associated with the Jordan algebra $V$. Let  $G(\Omega)$ be the subgroup of linear transformations of $V$ which preserve $\Omega$. Then $G(\Omega)$ is a reductive group, which acts transitively on $\Omega$. The same properties hold for its neutral component, which we denote by  $G_0$. The stabilizer $K_0=(G_0)_e$ of the unit $e$ is a maximal compact subgroup of $G_0$ and it is the neutral component of  $Aut(V)$, the {\it automorphism group} of the Jordan algebra $V$. Moreover $K_0=G_0\cap O(V)$, where $O(V)$ is the orthogonal group of the inner product on $V$. The space $\Omega\simeq G_0/K_0$ is a Riemannian symmetric space.\\

Let $c$ be an {\it idempotent} element\index{Idempotent element} in $V$ : $c^2=c$. Then the only possible eigenvalues of $L(c)$ are $1$, $\frac 12$, $0$ and $V$ is the direct sum of the corresponding eigenspaces $V(c,1)$, $V(c,\frac 12)$ and $V(c,0)$. The decomposition
$$V=V(c,1)\oplus V(c,\frac 12)\oplus V(c,0),$$
is called the {\it Peirce decomposition\/} \index{Decomposition!Peirce} of $V$ with respect to the idempotent $c$.  It is an orthogonal decomposition with respect to any scalar product satisfying (\ref{assoc}).
 
Two idempotents $c$ and $d$ are said to be {\it orthogonal\/} if $(c|d)=0$, which is equivalent to $cd=0$. An idempotent is said to be {\it primitive\/} if it is not the sum of two non-zero idempotents. An idempotent $c$ is primitive if and only if $\dim V(c,1)=1$. We say that $(c_j)_{1\leq j\leq m}$ is a {\it Jordan frame}\index{Jordan!frame} if each $c_j$ is a primitive idempotent and
$$c_ic_j=0,\;\; i\not= j$$
$$c_1+c_2+\ldots+c_m=e.$$
All the Jordan frames have the same number of elements which we  denote by $r$. The integer $r$ is the {\it rank\/} of the Jordan algebra $V$.\\
The group $K$ acts transitively on the set of primitive idempotents, and also on the set of Jordan frames. Therefore if we fix a Jordan frame $(c_j)_{j=1}^r$, then every element $x\in V$ can be written in the form
$$x=k(\sum_{j=1}^r\lambda_jc_j)$$
where $k\in K$ and $\lambda_1,\ldots, \lambda_r$ real numbers. The scalars $(\lambda_j)_{1\leq j\leq r}$ are unique and called the {\it spectral values} \index{Spectral!values} of $x$. We define the {\it determinant}\index{Determinant} and the {\it trace} \index{Trace} of the Jordan algebra by
$$\det(x)=\prod_{j=1}^r\lambda_j,\;\;  \tr(x)=\sum_{j=1}^r\lambda_j.$$
 The trace is a linear form of $V$ and the determinant is a homogeneous polynomial on $V$ of degree $r$. One can show that
 $$V^\times=\{x\in V\; , \; \det(x)\not=0\}.$$
 From now we assume that the scalar product of $V$ is given by
\begin{equation}\label{inn_prod}
(x|y)=\tr(xy).
\end{equation}
 \begin{example}
  The vector space $V=Sym(r,\mathbb{R})$ of $r\times r$ real symmetric matrices is a Euclidean Jordan algebra with the product $x\circ y=\frac12(xy+yx)$ and the scalar product $(x|y)=\Tr(xy)$. The quadratic representation is given by $P(x)y=xyx$.  In this case the determinant and the trace are  the usual matrix determinant and trace. The corresponding symmetric cone is $\Omega=Sym^{++}(r,\mathbb{R})$ the set of definite positive symmetric matrices.  An idempotent is an orthogonal projection $c=\left(\begin{smallmatrix} I_p&0\\0&0\end{smallmatrix}\right)$ with $r=p+q$. Then
$$
V(c,1)=\{\left(\begin{smallmatrix} a&0\\0&0\end{smallmatrix}\right) \;;\; a \;:\; \textrm{$p\times p$ symmetric matrix}\},$$
$$V(c,\frac 12)=\{\left(\begin{smallmatrix} 0&d\\d^t&0\end{smallmatrix}\right) \;;\; d \;:\; \text{$p\times q$ matrix}\},$$
$$V(c,0)=\{\left(\begin{smallmatrix} 0&0\\0&b\end{smallmatrix}\right) \;;\; b \;:\; \textrm{$q\times p$ symmetric matrix}\}.$$
\end{example}

\newpage
 \section{Hermitian symmetric spaces of tube type}\label{section3}

{\it In this section we characterize irreducible Hermitian symmetric spaces
of tube type using the Jordan theory and give the classification of such spaces.}\\

By complexification, we get a complex Jordan  algebra $V_\mathbb{C}$. We denote the complex conjugation with respect to $V$ by $\eta$. We also use the notation
\begin{equation}\label{conjug}
\eta(z)=\bar z.
\end{equation}
We extend the inner product (\ref{inn_prod}) of $V$ to 
  the Hermitian inner product of $V_\mathbb{C}$ defined by
\begin{equation}\label{herm_inn}
(z|w)=\tr(z\bar w).
\end{equation}
Let $$S=\{\sigma\in V_\mathbb{C}\; ,\; \sigma^{-1}=\bar\sigma\}.$$
This is a connected compact sub-manifold of $V_\mathbb{C}$.
\begin{proposition}[{\cite{Faraut-Koranyi}}] 
For $\sigma\in V_\mathbb{C}$ the following properties are equivalent :
\begin{enumerate}
\item[$(i)$] $\sigma\in S$,
\item[$(ii)$] $\sigma=\exp(iu)$, where $u\in V$,
\item[$(iii)$] There exists a Jordan frame $(c_j)_{1\leq j\leq r}$ of $V$ and complex numbers $(\xi_j)_{1\leq j\leq r}$ of modulus 1, such that $\sigma=\sum_{j=1}^r\xi_jc_j$.
\end{enumerate}
\end{proposition}

Denote by  $\mathbb{L}=Str(V_\mathbb{C})$ be the {\it structure group} \index{Group!structure} of $V_\mathbb{C}$,  i.e. the set of $g\in GL(V_\mathbb{C})$ such that
$$ P(gz)=gP(z)g'$$
or equivalently
$$g(z\Box w)g^{-1}=(g z)\Box(g'^{-1}w).$$ 
Consider the group
$$L(S)=\{g\in GL(V_\mathbb{C})\; ,\; g(S)=S\}.$$
Then $L(S)=\mathbb{L}\cap U(V_\mathbb{C})$, where $U(V_\mathbb{C})$ is the unitary group of the Hermitian inner product (\ref{herm_inn}) on $V_\mathbb{C}$. It acts transitively on $S$. Moreover, the stabilizer of $e$ in $L(S)$ coincides with $Aut(V)$ (we extend the automorphisms of $V$ as complex linear automorphisms of $V_\mathbb{C}$).  The involution (\ref{conjug}) preserves $S$ and $e$ is its unique isolated fixed point. The set of fixed points of the corresponding involution of $L(S)$, $g\mapsto \eta\circ g\circ \eta$, is $L(S)_e=Aut(V)$. Hence, with the  metric induced by the Hermitian product (\ref{herm_inn}), $S$ is a Riemannian symmetric space of compact type isomorphic to $L(S)/Aut(V)$.\\
Let $U$ be the identity component of $L(S)$, and $U_e$ the stabilizer of $e$ in $U$. Then we have $Aut(V)^\circ\subset U\subset Aut(V)$, and
$$S\simeq U/U_e.$$
Let $(c_j)_{1\leq j\leq r}$ be a Jordan frame of $V$. Then every element $z\in V_\mathbb{C}$ can be written as in the form
$$z=u(\sum_{j=1}^r\lambda_jc_j),$$
where $u\in U$ and $0\leq \lambda_1,\ldots,\lambda_r$. The {\it spectral norm} \index{Specral!norm} of $z$ is then defined by
$$|z|=\sup_{1\leq j\leq r}\lambda_j.$$
It turns out to be a norm on $V_\mathbb{C}$, invariant under the group $U$. \\
Introduce the domain $D$ in $V_\mathbb{C}$ as th open unit ball for the spectral norm
$$D=\{z\in V_\mathbb{C}\; ,\; |z|<1\}.$$
Recall that the {\it Shilov boundary} \index{Shilov boundary} of $D$ is the smallest closed set in $\overline{D}$ where the principle of the maximum holds.

\begin{theorem}[{\cite{Faraut-Koranyi}}]
$D$ is a bounded symmetric domain \index{Boundex symmetric domain} and $S$ is its Shilov boundary.
\end{theorem}
There is a realization of the domain $D$ as a    tube domain through the {\it Cayley transform}.\index{Cayley transform}  Let $T_\Omega$ be the tube over the symmetric cone $\Omega$,
$$T_\Omega=V+i\Omega=\{z=x+iy\in V_\mathbb{C}\;,\; y\in\Omega\}.$$
The Cayley transform $c$ and its inverse $p$ are given (in their domains of definition) by
$$\begin{array}{rl}
p(z)&=(z-ie)(z+ie)^{-1}\\
c(w)&=i(e+w)(e-w)^{-1}.
\end{array}$$
\begin{proposition}[{\cite{Faraut-Koranyi}}]
The map $p$ induces a biholomorphic isomorphism from $T_\Omega$ onto $D$, and
\begin{equation*}
p(V)=\{\sigma \in S \; ,\; \det(e-\sigma)\not=0\}.
\end{equation*}
\end{proposition}
Both domains $T_\Omega$ and $D$ are biholomorphically equivalent, and $V$
can be thought of as the Shilov boundary of $T_\Omega$ and its image under
the transformation $p$ is an open dense in $S$.\\

Let $G=G(D)$ be the neutral component of the group of biholomorphic
diffeomorphisms of $D$. It is a semi-simple Lie group and the stabilizer of
$0\in D$ in $G$ is a maximal compact subgroup of $G$ which coincides with
$U$.\\

To describe the group $G$, we use the Cayley transform\index{Cayley transform}. Let $G^c:=G(T_\Omega)$ be the neutral component of the group of biholomorphic diffeomorphisms of $T_\Omega$. Then 
$$c^{-1}\circ G\circ c=G^c.$$
We already know some subgroups of $G^c$. In fact , an element of $G_0$ acts
on $T_\Omega$ and we can identify $G_0$ with a subgroup of $G^c$. \\

For $v\in V$, the translation
\begin{equation*}
t_v : z\mapsto z+v
\end{equation*}
is a holomorphic automorphism of $T_\Omega$ and the group of all real translations $t_v$ is an Abelian subgroup $N^+$ of $G^c$ isomorphic to the vector space $V$.\\
The {\it inversion}
\begin{equation*}
j : z\mapsto z^{-1}
\end{equation*}
belongs to $G^c$. We set $N^-=j\circ N^+\circ j$. It is the subgroup of $G^c$ of the maps
$$\widetilde{t}_v =j\circ t_v\circ j : z\mapsto (z^{-1}-v)^{-1}, \; \; v\in V,$$
and it is an Abelian subgroup of $G^c$ isomorphic to $V$.
\begin{theorem}[{\cite{Faraut-Koranyi}}]
The subgroups $G_0$ and $N^+$, together with the inversion $j$, generate $G^c$.
\end{theorem}
The semi-direct product $P^+=G_0N^+$ is a maximal parabolic subgroup of $G^c$. The homogeneous space $G^c/P^+$ is then a (real) compact manifold which contains $V$ as an open dense subset,
\begin{equation}\label{compactification}
V \to G^c/P^+ \, :\, v\mapsto g_vP^+,
\end{equation}
where $g_v(z)=j(z)+v$. The manifold $G^c/P^+$ is the {\it conformal compactification} \index{Conformal compactification} of the Jordan algebra $V$, and it is isomorphic to the Shilov boundary $S$ of $D$.

 \begin{example} If $V$ is the Jordan algebra $Sym(r,\mathbb{R})$,  then $D$ is the Siegel domain
 $$D=\{z\in Sym(r,\mathbb{C}),\; I_r-zz^*>>0\}.$$
 It is holomorphically isomorphic to the upper half domain
 $$T_\Omega=\{z=x+iy\in Sym(r,\mathbb{C}),\;  y\in Sym^{++}(r,\mathbb{R})\}.$$
 In this case, $G=Sp(r,\mathbb{R})/\{\pm Id\}$ , where $Sp(r,\mathbb{R})$ is the symplectic group.
 \end{example}

 
 Here we give the classification of tube domains, their Shilov boundaries and the corresponding Euclidean Jordan algebras.
 
 \begin{center}
\begin{table}[h]
\begin{tabular}[t]{cccc}
\hline
$V$  & $V_\mathbb{C}$ & $D\simeq G/U$ & $S\simeq U/U_e$ \\
\hline
$Sym(m,\mathbb{R})$ & $Sym(m,\mathbb{C})$ & $Sp(2m,\mathbb{R})/U(m)$ & ${U(m)/O(m)}$ \\
${Herm(m,\mathbb{C})}$ & ${Mat(m,\mathbb{C})}$ & ${SU(m,m)/S(U(m)\times U(m))}$& ${U(m)}$\\
${Herm(m,\mathbb{H})}$ & ${Skew(2m,\mathbb{C})}$ & ${SO^*(4m)/U(2m)}$ & ${U(2m)/SU(m,\mathbb{H})}$\\
${\mathbb{R}\times\mathbb{R}^{q-1}}$ & ${\mathbb{C}\times\mathbb{C}^{q-1}}$ &
${SO_0(2,q)/SO(2)\times SO(q)}$ & ${(U(1)\times S^{q-1})/\mathbb{Z}_2}$\\
${Herm(3,\mathbb{O})}$ & ${Mat(3,\mathbb{O})}$ & ${E_{7(-25)}/U(1)E_6}$ &
${U(1)E_6/F_4}$\\
\hline
\end{tabular}
\caption{Tube domains and their Shilov boundaries }
\end{table}
\end{center}
\begin{table}[h]
\begin{center}
\begin{tabular}[t]{cccc}
\hline
${V}$  & ${n}$ & ${r}$ & ${d}$\\
\hline
${Sym(m,\mathbb{R})}$                 & ${\frac{1}{2}m(m+1)}$ & ${m}$ & ${1}$\\
${Herm(m,\mathbb{C})}$               & ${m^2}$               & ${m}$ & ${2}$  \\
${Herm(m,\mathbb{H})}$               & ${m(2m-1)}$           & ${m}$ & ${4}$\\
${\mathbb{R}\times\mathbb{R}^{q-1}}$ & ${q}$                 & ${2}$   & ${q-2}$\\
${Herm(3,\mathbb{O})}$               & ${27}$                  & ${3}$   &${8}$\\
\hline
\end{tabular}
\caption{The dimension, rank and the Peirce invariant}
\end{center}
\end{table}

\newpage

\section{Cayley type symmetric spaces}\label{section4}\index{Symmetric space!Cayley type}
 
{\it In this section we characterize the causal symmetric spaces of
  Cayley type; we prove in particular that if $G/K$ is a Hermitian
  symmetric space of tube type, then $G/H$ is a causal symmetric space of
  Cayley type. }\\

 Let $\frak g$ be the Lie algebra of $G=G(D)$ and $\frak g^c$ be the Lie algebra of $G^c=G(T_\Omega)$.   Let $\frak g_0$ be the Lie $G_0$. The Lie algebra of $K_0$ is set of all derivations $Der(V)$ of $V$. Let 
 $\frak p_0=\{L(v) \; ,\; v\in V\}.$
 Then the Cartan decomposition is given by $\frak g_0=\frak k_0\oplus \frak p_0$. The Lie algebra of $U$ is $\frak u=\frak k_0\oplus i\frak p_0$. Let  $G_\mathbb{C}$ be the Lie group generated by $Str(V_\mathbb{C})$, the complex translation of $V_\mathbb{C}$ and  $j$. Then $G$ and $G^c$ are two real forms of $G_\mathbb{C}$.  The Lie algebra $\frak g_\mathbb{C}$ of $G_\mathbb{C}$ is the set of vector fields $X$ on $V_\mathbb{C}$  of the form
 $$X(z)=u+Tz-P(z)v,$$
 with $u, v\in V_\mathbb{C}$ and $T\in \frak{str}(V_\mathbb{C})$, where $\frak{str}(V_\mathbb{C})$ is the Lie algebra of the structure group of $V_\mathbb{C}$ and coincides with $(\frak g_0)_\mathbb{C}$. If
  \begin{eqnarray*}
 X_1(z)&=&u_1+T_1z-P(z)v_1,\\
 X_1(z)&=&u_2+T_2z-P(z)v_2,
 \end{eqnarray*} 
 then the bracket $[X_1,X_2]$ is given by
 $$[X_1,X_2](z)=u+Tz-P(z)v,$$
 with
 \begin{eqnarray*}
 u&=&T_1u_2-T_2u_1,\\
 T&=&[T_1,T_2]+2(u_1\Box v_2)-2(u_2\Box v_1),\\
 v&=&-T^*_1v_2+T^*_2v_1.
 \end{eqnarray*} 
 A vector field $X$ in $\frak g_\mathbb{C}$ 
 $$X(z)=u+Tz-P(z)v,$$
can be identified with $(u,T,v)\in V_\mathbb{C}\times\frak{str}(V_\mathbb{C})\times V$, then
$$\frak g_\mathbb{C}\simeq V_\mathbb{C}\times \frak{str}(V_\mathbb{C})\times V_\mathbb{C}.$$
With this identification we have,
 \begin{eqnarray*}
 \frak g_0&=&\{(u,T,u) \; ,\; u\in V, T\in \frak k_0\}\simeq V\times \frak k_0\times V,\\
  \frak g^c&=&\{(u,T,v) \; ,\; u,\, v\in V, T\in \frak g_0\}\simeq V\times \frak g_0\times V,\\
 \frak g&=&\{(w,T,\bar w) \; ,\; w\in V_\mathbb{C}, T\in \frak u\},\\
 \frak u&=&\{(u,T,-u) \; ,\; u\in V, T\in \frak k_0\}.
\end{eqnarray*} 
Define $X_0=(0,I,0)$, then $\textrm{ad}(X_0)$ has eigenvalues $1, 0, -1$ with the eigenspaces
\begin{eqnarray*}
\frak g^c_1&=&\{(u,0,0) \; ,\; u\in V\}\simeq V,\\
\frak g^c_0&=&\{(0,T,0) \; ,\; T\in\frak g_0\}\simeq \frak g_0,\\
\frak g^c_{-1}&=&\{(0,0,v) \; ,\; v\in V\}\simeq V.\\
 \end{eqnarray*} 
 The decomposition
 $$\frak g^c=\frak g^c_1+\frak g^c_0+\frak g^c_{-1}$$
 is called the {\it Kantor-Koecher-Tits}\index{Decomposition!Kantor-Koecher-Tits} decomposition of $\frak g^c$.
     
  Consider the involutions $\sigma^c$ and $\theta^c$ of $G^c$ given by
 \begin{eqnarray*}
 \sigma^c(g)&=&(-j)\circ g\circ (-j)\\
 \theta^c(g)&=&j\circ g\circ j.
 \end{eqnarray*}
 We keep the same notation for the corresponding involutions of the Lie algebra $\frak g^c$. 
 \begin{proposition}[\cite{Koufany}, \cite{Koufany1}, \cite{Koufany2}] $\theta^c$ is the Cartan involution of $G^c$. It commutes with the involution $\sigma^c$. If $X=(u,T,v)\in \frak g^c$, then
 \begin{eqnarray*}
 \sigma^c(X)&=&(v,-T^*,u),\\
 \theta^c(X)&=&(-u,-T^*,-v).
 \end{eqnarray*}
 \end{proposition}
 Similarly, consider the involution $\sigma$ and $\theta$ of $G$ given by
 \begin{eqnarray*}
 \sigma(g)&=&\nu\circ g\circ \nu\\
 \theta^c(g)&=&(-\nu)\circ g\circ (-\nu).
 \end{eqnarray*}
 where $\nu(z)=\bar z$, and keep the same notation for the corresponding involutions of the Lie algebra $\frak g$. 
\begin{proposition}[\cite{Koufany}, \cite{Koufany1}, \cite{Koufany2}] $\theta$ is the Cartan involution \index{Cartan involution} of $G$. It commutes with the involution $\sigma$. If $X=(w,T,\bar w)\in \frak g$, then
 \begin{eqnarray*}
 \sigma(X)&=&(\bar w,\bar T,w),\\
 \theta(X)&=&(-w,-T,-\bar w).
 \end{eqnarray*}
 \end{proposition}
Let $\tau$ be the element of $G^c$ given by
$$\tau(z)=(e+z)(e-z)^{-1}.$$
Then its inverse is
$$\tau^{-1}(z)=(z-e)(z+e)^{-1}.$$
Notice that $c(z)=i\tau(z)$ and $c^{-1}(z)=\tau^{-1}(iz)$.  For any element $g\in G_0$,
$$c^{-1}\circ g\circ c=\tau \circ (g^*)^{-1}\circ \tau^{-1}.$$
Since $G_0$ is a reductive group, we have
$$c^{-1}\circ G_0\circ c=\tau \circ G_0\circ \tau^{-1}.$$
Let $H=G^{\sigma}$  and $H^c=(G^c)^{\sigma^c}$.  Then we have
\begin{theorem}[\cite{Koufany}, \cite{Koufany1}, \cite{Koufany2}]
\begin{enumerate}
\item $H=c^{-1}\circ G_0\circ c$ and $H^c=\tau \circ G_0\circ \tau^{-1}$.
\item The subgroups $H$ and $H^c$ coincides and $H=H^c=G\cap G^c$.
\item The symmetric space $\mathcal{M}=G/H\simeq G^c/H^c$ is a Cayley type symmetric space and any symmetric space of Cayley type is given in this way.
\end{enumerate}
\end{theorem}
 More precisely, we have
 \begin{eqnarray*}
\frak h^c&=&\{(u,T,u) \; ;\; u\in V, T\in \frak k_0\},\\
\frak q^c&=&\{(u,L(v),-u) \; ;\; u, v\in V\},\\
\frak k^c&=&\{(u,T,-u) \; ;\; u\in V, T\in \frak k_0\},\\
\frak p^c&=&\{(u,L(v),u) \; ;\; u, v\in V\}.
\end{eqnarray*}
Let $C_1$, respectively $C_2$,  be the cone in $\frak q^c$  given by
$$C_1=\{(u,2L(v),-u) \; ;\; (u+v)\in-\overline{\Omega}, (u-v)\in\overline{\Omega}\},$$
respectively
$$C_2=\{(u,2L(v),-u) \; ;\; (u+v)\in\overline{\Omega}, (u-v)\in\overline{\Omega}\}.$$
$C_1$ and $C_2$ are two $\text{Ad}(H^c)-$invariant, regular cones isomorphic to $\overline{\Omega}\times\overline{\Omega}$.
Moreover
$$C_1\cap \frak p^c\not=\emptyset, \; C_1\cap \frak k^c=\{0\}$$
and
$$C_2\cap \frak k^c\not=\emptyset, \; C_1\cap \frak p^c=\{0\}.$$
Thus $C_1$ (respectively $C_2$) defines a non-compactly (resp. compactly) causal structure on $G^c/H^c$.

 \newpage
  
  \section[The 2-transitivity and Cayley type spaces]{The 2-transitivity
    property on $S^2_\top$ and Cayley type symmetric
    spaces}\label{section5}\index{Symmetric space!Cayley type}

{\it  In this section we will give a causal compactification of cuasal
  symmetric spaces of Cayley type.}\\

  \begin{definition}
  Let $\mathcal{M}$ be a causal $G-$manifold. A causal compactifiction \index{Causal compactification} of $\mathcal{M}$ is a pair $(\mathcal{N},\Phi)$ such that\begin{enumerate}
\item $\mathcal{N}$ is a compact causal $G-$manifold.
\item The map $\phi : \mathcal{M}\to \mathcal{N}$ is causal.
\item The map $\Phi$ is $G-$equivariant, i.e., $\Phi(g\cdot x)=g\cdot\Phi(x)$, for every $g\in G$ and every $x\in \mathcal{M}$.
\item $\Phi(\mathcal{M})$ is open and dense in $\mathcal{N}$.
\end{enumerate}

  \end{definition}
  
  Two points $z, w\in V_\mathbb{C}$ are called {\it transversal},\index{Transversal} and we write $z\top w$, if and only if $\det(z-w)\not=0$. This condition is equivalent to $\Det P(z-w)\not=0$.  We denote the set of transversal elements in $S^2$ by
 $$S^2_\top=\{(\sigma,\tau)\in S^2\; ,\; \sigma\top\tau\}=\{(\sigma,\tau)\in S^2\; ,\; \det(\sigma-\tau)\not=0\},$$
   
  
  \begin{theorem}[\cite{Koufany}, \cite{Koufany1}, \cite{Koufany2}] The group $G$ acts transitively on $S^2_\top$. The stabilizer of the element $(e,-e)\in S^2_\top$ in $G$ is the group $H=c^{-1}\circ G_0\circ c$.
\end{theorem}
Hence, the Cayley symmetric space is $G-$equivariant to $S^2_\top$,
$$G/H\simeq S^2_\top.$$
Since  $S^2_\top$ is open dense in $S^2$, and since $S^2$ is a compact causal $G-$manifold (because we already prove that $S$ is a causal $G-$manifold), the manifold $S^2$ is a causal compactification of $G/H$.\\
The "non bounded" realization of $\mathcal{M}=G^c/H^c$ is such that
$$\mathcal{M}\cap(V\times V)=\{(x,y)\in V\times V \; ;\; \det(x-y)\not=0\}.$$

\begin{examples}
\begin{enumerate}
\item If $\mathcal{M}=SU(n,n)/GL(n,\mathbb{C})\mathbb{R}^+$. Then 
$$D=SU(n,n)/S(U(n)\times U(n))=\{z\in Mat(n,\mathbb{C}) \; ; \; I_n-z^*z\gg 0\},$$
its Shilov boundary is $S=U(n)$ and
$$\mathcal{M}\simeq\{(z,w)\in U(n)\times U(n) \; ;\; \Det(z-w)\not=0\}.$$
\item If $\mathcal{M}=Sp(n,\mathbb{R})/GL(n,\mathbb{R})\mathbb{R}^+$. Then
$$D=Sp(n,\mathbb{R})/U(n)=\{z\in Sym(n,\mathbb{C}) \; ; \; I_n-z^*z\gg 0\},$$
its Shilov boundary is the Lagrange Grassmann manifold $S=U(n)/O(n)=\{z\in U(n) \; ;\; z^t=z\}$ and
$$\mathcal{M}\simeq \{(z,w)\in U(n)\times U(n) \; ;\; z^t=z, w^t=w, \Det(z-w)\not=0\}.$$
\end{enumerate}
\end{examples}

\newpage

\section[The causal Lie semigroup]{The Lie semigroup associated with the Cayley type symmetric space}\label{section6}

{\it In this section we inverstigate the semigroup $S_\Omega$ of compressions of
  $\Omega$. We prove in particular a triple decomposition and that
  $S_\Omega$ is the real part of the
  holomorphic semigroup of compression of the tube domain. We give a new
  characterization of the Riemannian metric of $\Omega$ and prove that
  $S_\Omega$ is a semigroup of contractions of this metric. The
  Hilbert metric on $\Omega$ is also studied}

\subsection{The compression semigroup of the Hermitian domain $D$}\index{Semi-group!of compressions}
  Let  $C^c_{\textrm{max}}$ be the maximal cone \index{Cone!maximal} in $\frak g^c$. It is the closed convex cone given by
 $$C^c_{\textrm{max}}=\{X\in \frak g^c \; ;\; X(v)\in\overline{\Omega}, \forall v\in V\}.$$
 The cone $C_{\textrm{max}}=\text{Ad}(c)(C^c_{\textrm{max}})$ is the maximal cone in $\frak g$.
 Let 
 \begin{equation*}
 \Gamma(C_{\textrm{max}})=G\exp(iC_{\textrm{max}}).
 \end{equation*}
 The following theorem is due to Olshanski{\u \i}
 \begin{theorem}[\cite{Ol2}]
  The set $\Gamma(C_{\textrm{max}})$ is a Lie semigroup \index{Semi-group!Lie} (associated with $C_{\textrm{max}}$) and it is the semigroup of compressions of the Hermitian domain $D$,
 $$\Gamma(C_{\textrm{max}})=\{g\in G_\mathbb{C} \; ;\; g(D)\subset D\}.$$
 Moreover
 $$\Gamma(C_{\textrm{max}})^\circ=\{g\in G_\mathbb{C} \; ;\; g(\overline{D})\subset D\}.\\$$
   \end{theorem}
\
The convex cone $C_1$ is the maximal cone $c^c_{\textrm{max}}$ in $\frak q^c$. It defines a non-compactly causal structure on $\mathcal{M}$. Therefore $\mathcal{M}$ is and ordering symmetric space 
Let $\Gamma$ the semigroup associated with the order of $\mathcal{M}$,
$$S_\succeq=\{g\in G^c \; ;\; g(e,-e)\succeq (e,-e)\}.$$ Then we have
\begin{theorem}[\cite{Koufany}, \cite{Koufany1}, \cite{Koufany2}]
The semigroup $S_\succeq$ satisfies 
\begin{enumerate}
\item $S_\succeq=\exp(c^c_{\textrm{max}})H$.
\item $\Gamma(C_{\textrm{max}})\cap G^c=S_\succeq^{-1}$
\end{enumerate}
\end{theorem}

\subsection{The compression semigroup of the symmetric cone $\Omega$}
Recall that when $P$ is the parabolic subgroup $P=G_0N^+$, the compact symmetric space $\mathcal{X}=G^c/P$ is conformal compactification of $V$, and the imbedding of $V$ into $\mathcal{X}$ given by (\ref{compactification}) is  an open dense embedding. A Lie semigroup that is naturally related to the action of $G^c$ on $\mathcal{X}$ occurs as the semigroup of compressions of $\Omega$ in $G^c$ :
\begin{equation}
S_\Omega=\{\gamma\in G^c \; ;\; \gamma \Omega\subset \Omega\}.
\end{equation}
Since the closure $\widetilde{\Omega}$ of $\Omega$ in $\mathcal{X}$ is compact with $\Omega$ as interior, the compression semigroup $S_\Omega$ is a closed semigroup of $G^c$. Moreover $S_\Omega$ contains $G_0$, and its interior is
$$S_\Omega^\circ=\{\gamma\in G^c \; ;\; \gamma \widetilde{\Omega}\subset\Omega\}.$$

Now let
\begin{eqnarray*}
S_\Omega^+&=&\{\gamma^+_v : z\mapsto z+v \; ;\; v\in\overline{\Omega}\},\\
S_\Omega^-&=&\{\gamma^-_v : z\mapsto (z^{-1}+v)^{-1}\; ;\; v\in\overline{\Omega}\}.
\end{eqnarray*}
Then it is easy to see that $S_\Omega^\pm$ and $G_0$ are closed sub-semigroups in $S_\Omega$. Hence $S_\Omega^+G_0 S_\Omega^-\subset S_\Omega$.

\begin{theorem}[\cite{Koufany}, \cite{Koufany2}]
\begin{enumerate}
\item The compression semigroup $S_\Omega$ is equal to the semigroup $S_\succeq$
\item The sub-semigroups $S^+_\Omega$ and $S^-_\Omega$, together with the subgroup $G_0$, generate $S_\Omega$. More precisely, one has the following decomposition
 \begin{equation}
S_\Omega=S^+_\Omega G_0S^-_\Omega=N^+G_0N^-\cap S_\Omega.
\end{equation}
\end{enumerate}
\end{theorem}

If $\gamma=\gamma^+_ug\gamma^-_v\in S_\Omega$, then we write
\begin{equation}\label{not_cont}
n^+(\gamma):=u,\; A(\gamma):=g\;\mathrm{and},\; n^-(\gamma):=v.
 \end{equation}

\subsection{The semigroup of contractions}\index{Semi-group!of contractions}
The family of bilinear forms $\mathrm{g}_x$ given by,
\begin{equation*}
\mathrm{g}_x(u,v)=(P(x)^{-1}u|v), \quad x\in \Omega, u, v\in V,
\end{equation*}
  defines a $G(\Omega)-$invariant Riemannian metric \index{Riemannian metric} on $\Omega$, see \cite[Theorem III.5.3]{Faraut-Koranyi}. Therefore, $\Omega$ is a 
Riemannian symmetric space isomorphic to $G(\Omega)_{\!\circ}/K(\Omega)_\circ$.

\begin{theorem}[\cite{Koufany3}]\label{distance_theorem}
Let $x, y\in\Omega$. Then there exists a unique curve of shortest length 
joining $x$ and $y$. The length of this curve is given by
\begin{equation*}
\delta(x,y)=\Bigl(\sum_{k=1}^r\log^2\lambda_k(x,y)\Bigr)^{1/2},
\end{equation*}
where $\lambda_1(x,y),\ldots,\lambda_r(x,y)$ are the spectral values of 
of $P(y)^{-1/2}x$.
\end{theorem}
$\delta(x,y)$ is the Riemannian distance \index{Riemannian distance} of $x$ and $y$, and the scalars
$$\mu_k(x,y):=\log^2(\lambda_k(x,y)$$
are by definition the {\it angles} \index{Angles} or the {\it compounds distance}.\index{Compounds distance}\\
Using the notations (\ref{not_cont}), we set
\begin{equation*}
S_1=\{\gamma\in S_\Omega \tq n^+(\gamma)\in\Omega\},
\end{equation*}
and
\begin{equation*}
S_2=\{\gamma\in S_\Omega \tq n^-(\gamma)\in\Omega\}. 
\end{equation*}

\begin{theorem}[\cite{Koufany3}]\label{contractio_theoerem}
Let $k\in\{1,\ldots,r\}$. The following holds:
\begin{enumerate}
\item  For any $\gamma\in S$ and for any $x, y\in\Omega$ : $\mu_k(\gamma\cdot x,\gamma\cdot y)\leq\mu_k(x,y).$
\item For any $\gamma\in S_1\cup S_2$ and for any $x, y\in\Omega$ : $\mu_k(\gamma\cdot x,\gamma\cdot y)<\mu_k(x,y).$
\item For any $\gamma\in S_1\cap S_2$, there exists $\kappa(\gamma)$, $0<\kappa(\gamma)<1$, such that for any $x, y\in\Omega$ : $\mu_k(\gamma\cdot x,\gamma\cdot y)\leq\kappa(\gamma)\mu_k(x,y).$
\end{enumerate}
\end{theorem}

 As an easy consequence, Theorem \ref{contractio_theoerem} implies that the elements of the semigroup $S_\Omega$ are contractions of the distance $\delta$. More precisely we have
\begin{corollary}[\cite{Koufany3}]
The following holds:
\begin{enumerate}
\item For any $\gamma\in S_\Omega$, and $x, y\in\Omega$ : $\delta(\gamma\cdot x,\gamma \cdot y)\leq\delta(x,y).$
\item For any $\gamma \in S_1\cup S_2$ and $x, y\in\Omega$ : $\delta(\gamma\cdot x,\gamma \cdot y)<\delta(x,y).$
\item For any $\gamma\in S_1\cap S_2$, there exists $\kappa(\gamma)$, 
$0<\kappa(\gamma)<1$, such that, for all $x, y\in\Omega$ : $\delta(\gamma\cdot x,\gamma \cdot y)\leq \kappa(\gamma)\;\delta(x,y).$
\end{enumerate}
\end{corollary}
  
\subsection{Hilbert's projective metric}
Let $E$ be a real Banach space and $C$ be a closed convex pointed cone,
where pointed means $C\cap -C=\{0\}$. The relation $\preccurlyeq$ is defined on $E$ by saying that $x\preccurlyeq y$ if and only if $y-x\in C$. \\
For $x\in E$ and $y\in \text{int}(C)$ we let
\begin{equation*}
M(x,y):=\inf\{\lambda \;|\; x\preccurlyeq \lambda y\},
\end{equation*}
and
\begin{equation*}
m(x,y):=\sup\{\mu \;|\; \mu y\preccurlyeq x\}.
\end{equation*}
Hilbert's projective metric\index{Hilbert's projective metric} is defined on $\text{int}(C)$ by
\begin{equation}\label{hilbert_metric}
d(x,y)=\log\frac{M(x,y)}{m(x,y)}.
\end{equation}

In the case of $\mathbb{R}^n_+$, Hilbert's projective metric is
$d(x,y)=\log\frac{\mathrm{max}{\frac{x_i}{y_i}}}{\mathrm{min}{\frac{x_i}{y_i}}}$
where $x=(x_1,\ldots,x_n)$ and $y=(y_1,\ldots,y_n)$ are two vectors of
$\mathbb{R}^n_+$.

The Hilbert projective metric may be applied to variety of problems
involving positive matrices and positive integral operators. For example
one can use it to solve some Volterra equations. It is also particularly
useful in proving the existence of the fixed point for positive
operators defined in a Banach space. In this way, it has been shown by
Bushell \cite{Bushell3} that Hilbert projective metric may be applied to
prove that, if $T$ a real nonsingular $r\times r$ matrix, then there exists
a unique real positive definite symmetric $r\times r$ matrix $A$ such that 
\begin{equation}\label{equa_bushell}
T'AT=A^2.
\end{equation}
Notice that if $T$ is neither symmetric nor orthogonal the existence and
the uniqueness of $A$  is not an elementary problem, even if $r=2$. 

We will formulate the Hilbert projective metric on symmetric cones in a way most convenient for our purpose using Jordan algebra theory
and extend Bushell's Theorem to this class of convex cones. 

If we consider the cone $\Omega_\mathrm{Sym}$ of real symmetric positive
definite $r\times r$ matrices, then one can easily express the Hilbert
projective metric (\ref{hilbert_metric}) in terms of eigenvalues of
elements of $\Omega_\mathrm{Sym}$. Indeed, if $A$ and $B$ are in $\Omega_\mathrm{Sym}$, then 
$$M(A,B):=\inf\{\lambda \;|\; \lambda B-A\ord 0\}=\max_{\|x\|=1}\frac{(Ax|x)}{(Bx|x)},$$
and 
$$m(A,B):=\sup\{\lambda \;|\; \lambda B-A\ord 0\}=\min_{\|x\|=1}\frac{(Ax|x)}{(Bx|x)},$$ which are respectively the
greatest and the least eigenvalue of $B^{-1}A$. Observe that eigenvalues of
the matrix $B^{-1}A$ are the same of the
matrix $B^{-\frac{1}{2}}AB^{-\frac{1}{2}}=P(B^{-\frac{1}{2}})A$.

More generally, for symmetric cones, Hilbert's projective metric can be also formulated in terms of extremal eigenvalues : let $x$ and $y$ be in $\Omega$ and let $\lambda_M(x,y)>0$ and $\lambda_m(x,y)>0$ denote the greatest and the least eigenvalue of the element $P(y^{-\frac{1}{2}})x\in \Omega$. Then one can prove that 
\begin{equation*}
\lambda_M(x,y)=\max_{c\in\mathcal{J}(V)}\frac{(x|c)}{(y|c)},
\end{equation*}
and
\begin{equation*}
\lambda_m(x,y)=\min_{c\in\mathcal{J}(V)}\frac{(x|c)}{(y|c)},
\end{equation*}
see  \cite[Thoerem 4.2]{Koufany3}. Consequently, we have :
\begin{proposition}[\cite{Koufany4}]\label{def_hil}
If $x,\, y\in\Omega$, then the Hilbert metric of $x$ and $y$ is given by
\begin{equation}\label{chara_hilb}
d(x,y)=\log\frac{\lambda_M(x,y)}{\lambda_m(x,y)}=\log\bigl[\lambda_M(x,y)\lambda_M(y,x)\bigr].
\end{equation}
\end{proposition}

Furthermore, one can prove (see \cite{Koufany4}) that $(\Omega,d)$ is a pseudo-metric space. In other words, for any $x,\, y,\, z\in\Omega$, the following holds,
\begin{enumerate}
\item[(a)] $d(x,y)\geq 0$
\item[(b)] $d(x,y)=d(y,x)$
\item[(c)] $d(x,z)\leq d(x,y)+d(y,z)$
\item[(d)] $d(x,y)=0\Leftrightarrow \exists \lambda>0 : x=\lambda y$.
\end{enumerate}
Now the characterization (\ref{chara_hilb}) of the Hilbert metric allows us to prove the completeness :
\begin{proposition}[\cite{Koufany4}]\label{complete}
$(\Omega\cap S(V),d)$ is a complete metric space.
\end{proposition}
As application, we prove a generalization of the Bushell theorem \index{Bushell theorem}:
\begin{theorem}[\cite{Koufany4}]
Let $g\in G(\Omega)$ and $p\in\mathbb{R}$ such that $|p|> 1$. Then there exists a unique element $a$ in $\Omega$ such that 
$g(a)=a^p$. 
\end{theorem}

\newpage

\section[The transversality index]{The 2-transitivity property on $S^2$ and
  the transversality index}\label{section7}\index{Index!transversality}

{\it We introduce here, the transversality index, a new invariant on the
  Shilov boundary which characterize the action of $G$ on $S\times S$. In the
  particular case of symmetric matrices this invariant has been studied by Hua.}\\

Fix a Jordan frame $(c_j)_{1\leq j\leq r}$ and for $k=0,1\ldots,r$ let
\begin{equation*}
\epsilon_0=-e,\; \epsilon_k=\sum_{j=1}^kc_j-\sum_{j=k+1}^rc_j,\; \epsilon_r=e.
\end{equation*}
\begin{proposition}
There are exactly $r+1$ orbits in $V^\times$ under the action of $G_0$. The elements $\epsilon_k$, $0\leq k\leq r$ are the set of representatives of all the orbits.
\end{proposition}
Using this proposition and the Cayley transform one can prove the following
\begin{theorem}
There are exactly $r+1$ orbits in $S\times S$ under the action of $G$ represented by the family $(e,\epsilon_k)$, $0\leq k\leq r$.
\end{theorem}

Let $(\sigma,\tau)\in S\times S$. The {\it transversality index} of the pair $(\sigma,\tau)$ is defined by
\begin{equation}\label{transv_index}
\mu(\sigma,\tau)=k
\end{equation}
where $k$ is the unique integer, $0\leq k\leq r$ such that $(\sigma,\tau)$ is conjugate under $G$ to the pair $(e,\epsilon_k)$.\\

The transversality index can also be understood as follows : Recall that the rank $\rank(x)$ of an element $x\in V$ is by the number of its non-zero spectral values with their multiplicities counted. This is an invariant under the action of $G_0$. Observe that $x$ and $y$ have the same rank if and only if $P(x)$ and $P(y)$ have the same rank.\\
 
 \begin{proposition}[\cite{Clerc-Koufany}]
 Let $\sigma, \tau \in S$, then there exists $u\in U$ such that $u(\sigma)$ and $u(\tau)$ are transversal to $e$. In addition the integer
 $$\rank[c(u(\sigma))-c(u(\tau))]$$
 does not depends on the element $u\in U$ and
$$\mu(\sigma, \tau)=r-\rank[c(u(\sigma))-c(u(\tau))].$$
 \end{proposition} 
 Notice that 
$$\mu(\sigma,\tau)=k \iff \rank\, P(\sigma-\tau)=k+\frac{k(k-1)}{2}d.$$
In particular, 
$$  \mu(\sigma,\tau)=0 \iff \sigma\top \tau .$$

\newpage

\section[The triple Maslov index]{The 3-transitivity property on $S^3_\top$ and the triple Maslov
  index}\label{section8}\index{Index!Maslov}

{\it In this section we recall another invariant which characterize the
  action of $G$ on $(S\times S\times S)^\top$. This invariant was
  introduced by Clerc and {\O}rsted and provides a generalization of the triple Malsov index.}\\

 Let $S^3_\top$ be set of pairwise transversal elements in  $S^3$,
 $$S^3_\top=\{(\sigma_1,\sigma_2,\sigma_3)\in S^3\; ,\; \sigma_i\top\sigma_j,\; 1\leq i\not=j\leq j\}.\}.$$
Choose a Jordan frame $(c_j)_{1\leq j\leq r}$, $\sum_{j=1}^rc_j=e$, then we have
\begin{theorem}[\cite{Clerc-Orsted-1}]
There are exactly $r+1$ orbits in $S^3_\top$ under the action of $G$, represented by the family $(e,-e,-i\epsilon_k)$, $0\leq k\leq r$.
\end{theorem}
Let $(\sigma_1,\sigma_2,\sigma_3)\in S^3\top$. The {\it triple Maslov index} $\imath(\sigma_1,\sigma_2,\sigma_3)$ of the triplet $(\sigma_1,\sigma_2,\sigma_3)$ is defined by
\begin{equation}\label{triple-Maslov}
\imath(\sigma_1,\sigma_2,\sigma_3)=k-(r-k)=2k-r
\end{equation}
where $k$ is the unique integer $0\leq k\leq r$ such that $(\sigma_1,\sigma_2,\sigma_3)$ is conjugate under $G$ to the triplet $(e,-e,-i\epsilon_k)$.
One can prove that the triple Maslov index $\imath$ satisfies :
\begin{proposition}[\cite{Clerc-Orsted-1}]  The triple Maslov index is 
\begin{enumerate}
\item an integer valued function : for all  $(\sigma_1,\sigma_2,\sigma_3)\in S^3\top$, 
$$-r\leq \imath(\sigma_1,\sigma_2,\sigma_3)\leq r$$
\item invariant under the action of $G$ : for all  $(\sigma_1,\sigma_2,\sigma_3)\in S^3\top$ and all $g\in G$, 
\begin{equation*}
\imath(g(\sigma_1),g(\sigma_2),g(\sigma_3))=\imath(\sigma_1,\sigma_2,\sigma_3)
\end{equation*}
\item skew symmetric : for all  $(\sigma_1,\sigma_2,\sigma_3)\in S^3\top$ and all  permutation $\pi$ of $\{1,2,3\}$,
\begin{equation*}
\imath(\sigma_{\pi(1)},\sigma_{\pi(2)},\sigma_{\pi(3)})=\imath(\sigma_1,\sigma_2,\sigma_3)
\end{equation*}
\item a cocycle : for all $\sigma_1,\sigma_2,\sigma_3,\sigma_4\in S$ such that $\sigma_i\top\sigma_j$, $1\leq i\not=j\leq 4$,
\begin{equation*}
\imath(\sigma_1,\sigma_2,\sigma_3)=\imath(\sigma_1,\sigma_2,\sigma_4)+\imath(\sigma_2,\sigma_3,\sigma_4)+\imath(\sigma_3,\sigma_1,\sigma_4).
\end{equation*}
\end{enumerate}
\end{proposition}
There is another construction of the triple Maslov index as the integral of the Kaehler form of the Hermitian domain $D$. Let $z_1, z_2, z_3\in D$. Form the oriented geodesic triangle \index{Geodesic triangle} $\Delta(z_1, z_2, z_3)$, and consider any surface $\Sigma$ in $D$ which has this triangle as boundary. Let $\omega$ be the Kaehler form of the domain $D$. Then the real number  
$$\varphi(z_1, z_2, z_3)=\int_\Sigma \omega$$
is not depending on $\Sigma$, since the Kaehler form is closed, and is called the {\it symplectic area} \index{Symplectic area} of the triangle $\Delta(z_1, z_2, z_3)$. As the Kaehler form is invariant under $G$, this gives and invariant for the oriented triples in $D$. Now for $(\sigma_1, \sigma_2, \sigma_3)\in S^3_\top$, then
the limit
$$\lim_{z_j\to\sigma_j}\varphi(z_1,z_2,z_3)$$
exists and is equals to the triple Maslov index $\imath(\sigma_1, \sigma_2, \sigma_3)$, see \cite{Clerc-Orsted-2}. This definition of the triple Maslov index extends for general $(\sigma_1, \sigma_2, \sigma_3)\in S^3$, without the transversality condition. This requires a notion of {\it radial convergence}, 
 see \cite{Clerc2} for more details.

\newpage

 \section[The universal covering of $S$]{The universal covering of the Shilov boundary} \label{section9}\index{Universal
   covering!of the Shilov Boundary}
{\it In this section we introduce the causal structure of the Shilov
  boundary $S$ and give an explicit construction of the universal covering of $S$.}
 \subsection{The causal structure of the Shilov boundary}\label{section_causal}\index{Structure!causal} 
Usually, one demands that the cones in the causal structures be closed. For our purpose we will suppose them open. \\
The Jordan algebra $V$ has a natural  causal structure modelled after the symmetric cone $\Omega$. It is simply given by the symmetric cone $\Omega$ viewed as a causal cone  in the tangent space $T_x(V)=V$ at any point $x\in V$.  It is clear that this causal structure is $G^c-$invariant .\\
Consider $-e$ as the base point of $S$. The tangent space $T_{-e}(S)$ can
be identified to $iV$. Moreover $p(0)=-e$ and the Cayley transform
$p$\index{Cayley transform} is well defined in a neighbourhood   of $0$. Its derivative at $0$ is given by
$$Dp(0)=-2i\textrm{Id}_{V_\mathbb{C}}.$$
This allows us to transfer the causal structure form $V$ to $S$.  We define the causal cone $C_{-e}$ to be
$$C_{-e}=Dc(-e)(\Omega)=-i\Omega.$$
This is an invariant cone by the stabilizer of $-e$ in $G$. We can then define a $G-$invariant causal structure on $S$ as follow : Let $\sigma\in S$, then there exists $g\in G$ such that $g(-e)=\sigma$, then we define the causal cone $C_\sigma$ to be
\begin{equation}
C_\sigma = Dg(-e)(C_{-e})=Dg(-e)(-i\Omega).
\end{equation}
 The family $(C_\sigma)_{\sigma\in S}$ is the unique  $G-$invariant causal
 structure of the Shilov boundary $S$ modelled after $\Omega$.

\subsection{The construction of the universal covering}
The Shilov boundary $S\simeq U/U_e$ is not a semi-simple symmetric space. To construct the universal covering of $S$, we prefer to deal with its semi-simple part. \\
Consider the set
\begin{equation*}
S_\mathtt{1}=\{\sigma\in S \tq \det(\sigma)=1\},
\end{equation*}
then $S_\mathtt{1}$  is a connected sub-manifold of $S$.
There exists a character $\chi$ of the structure group $Str(V_\mathbb{C})$ such that
$$\det(gz)=\chi(g)\det(z)$$
for all $g\in Str(V_\mathbb{C})$ and all $z\in V_\mathbb{C}$.  Let 
\begin{equation*}
U_\mathtt{1} = \{u\in U\mid \chi(u) = 1\} .
\end{equation*} 
Then $U_\mathtt{1}$ is a compact semi-simple group and acts  transitively on $S_\mathtt{1}$ and
\begin{equation*}
S_\mathtt{1}\simeq U_\mathtt{1}^\circ/(U_\mathtt{1}^\circ\cap Aut(J)),
\end{equation*}
where $U_\mathtt{1}^\circ$ is the neutral component of the group $U_\mathtt{1}$. Thus we have
 \begin{proposition}[{\cite{Clerc-Koufany}}]
$S_\mathtt{1}$ a semi-simple Riemannian symmetric space of compact type.
\end{proposition}
The Lie algebra of $U_\mathtt{1}$ is $\mathfrak{u}_\mathtt{1} = \mathfrak{k}\oplus i\mathfrak{p}_1$ where
$$\mathfrak{p}_1 = \{L(v)\tq v\in V, \tr(v) = 0\}.$$ 
Let
\begin{equation*}
 \mathfrak{a}_1 = \{L(a)\tq a=\sum_{j=1}^r a_j c_j,\; a_j\in
\mathbb{R},\;\tr(a) = 0\}.
\end{equation*}
Then $ \mathfrak{a}_1$ is Cartan subspace of $ \mathfrak{p}_1$. We consider now the fundamental lattice \index{Lattice!fondamental} $\Lambda_0$ of $\mathfrak{a}_1$,
$$\Lambda_0=\text{lattice generated by }\{2\pi\frac{A_{j,k}}{(
  A_{j,k}|A_{j,k} )} \tq 1\leq j\not= k\leq r\},$$
where $A_{j,k}$ is the covector of $\frac{1}{2}(a_j-a_k)$. The unitary lattice \index{Lattice!unitary} $\Lambda$ of $\mathfrak{a}_1$ is given by
\begin{equation*}
\Lambda=\{H\in \mathfrak{a}_1 \tq \exp(iH)e=e\}.
\end{equation*}
One can  prove (see \cite{Clerc-Koufany}) that $\Lambda_0=\Lambda$.  According \cite[Theorem 3.6]{Loos} (see also  \cite[Ch. VII. Theorem 8.4 et Theorem 9.1]{Helgason2}) we have
\begin{theorem}[\cite{Clerc-Koufany}]
The symmetric space  $S_\mathtt{1}$ is simply connected.
\end{theorem}

Following a classical method, we will realize the universal covering of $S$. let 
\begin{equation*}
\widetilde{S} = \{(\sigma,\theta)\in S\times \mathbb{R}\tq \det(\sigma)= e^{ir\theta}\},
\end{equation*}
with the topology induced by the topology of $S\times \mathbb{R}$.  
\begin{theorem}[\cite{Clerc-Koufany}]
$\widetilde{S}$ is the universal covering of $S$.
\end{theorem}
In fact we prove that map
\begin{equation*}
S_\mathtt{1}\times \mathbb{R}\longrightarrow \widetilde{S}\quad (\sigma, \theta)\longmapsto (e^{i\theta}\sigma, \theta)
\end{equation*}
is a homeomorphism and bijective.\\

Now we will describe a covering of the conformal group $G$ and give and explicit action of it on the universal covering of the Shilov boundary.  \\
For $g\in G$ and $z\in D$ we define
$$j(g,z)=\chi(Dg(z)).$$
This is an element of the structure group $Str(V_\mathbb{C})$. it is easy to see that $j(g,z)\not=0$. Since $D$ is simply connected, we can find a determination $\varphi_g$ of the argument of $j(g,\cdot)$, that is
$$\forall z\in \mathcal{D},\;\; e^{i\varphi_g(z)}=\frac{j(g,z)}{|j(g,z)|}.$$
Two such determinations differs by $2\pi k$.\\
Consider the following group
$$ \Gamma = \bigl\{\big(g,\varphi_g\big)\;|\; g\in G\bigr\} .$$
The multiplicative law being given by
$$ (g,\varphi(g,\cdot))(h,\psi(h,\cdot))=\big(gh, \varphi(g,h(\cdot))+\psi(h,\cdot)\big).$$
Observe that $\Gamma$ can be identified with the closed sub-set of $G\times \mathbb{R}$ given by
$$ \big\{(g,\theta)\in G\times \mathbb{R}\;\mid\; e^{i\theta}= j(g,0)\big\}.$$
Thus $\Gamma$ becomes a topological group.
\begin{proposition}[\cite{Clerc-Koufany}]
For $ (g,\varphi_g)\in \Gamma$
and
$(\sigma,\theta)\in \widetilde{S}$, set
$$(g,\varphi_g)\cdot
(\sigma,\theta)=\big(g(\sigma),\theta+\frac{1}{r}\varphi(g,\sigma)\big) .$$
Then this defines a continuous action of  $\Gamma$ on $\widetilde{S}$.
\end{proposition}
To prove this proposition we need the following formula (see \cite[Lemme 3.6]{Clerc-Koufany})
\begin{equation}\label{equation_det_j}
\det(g(\sigma)) = \frac{j(g,\sigma)}{|j(g,\sigma)|} \det(\sigma), \;\; for g\in G, \; \sigma\in S.
\end{equation}
This requires the  causal structure of $S$.

\newpage

 \section{The Souriau index}\label{section10}\index{Index!Souriau}

{\it In this section we construct a {\it primitive} $m$ of the Maslov cocycle. That is an integer valued function 
$$m : \widetilde{S}\times \widetilde{S} \to \mathbb{Z}$$
which is skew symmetric and satisfies the following cohomology property
\begin{equation*}
\imath(\sigma_1,\sigma_2,\sigma_3)=m(\widetilde{\sigma}_1,\widetilde{\sigma}_2)+m(\widetilde{\sigma}_2,\widetilde{\sigma}_3)+m(\widetilde{\sigma}_3,\widetilde{\sigma}_1).
\end{equation*}
}

Let $\sigma$ belongs to the set $S_\top(-e)$ of elements in $S$ transversal to $-e$. We can define the logarithm of $\sigma$ by
\begin{equation}\label{log}
\log \sigma=\int_{-\infty}^0\big((se-\sigma)^{-1}-(s-1)^{-1}e\big)ds\in V_\mathbb{C}.
\end{equation}
The function $\log$ has the following standard properties :
\begin{proposition}[\cite{Clerc-Koufany}]
For any $\sigma\in S_\top(-e)$, 
\begin{enumerate}
\item[$(i)$] $\exp(\log \sigma )=\sigma.$
\item[$(ii)$] $e^{\tr(\log \sigma )}=\det(\sigma).$
\item[$(iii)$] $\log \sigma^{-1} =-\log \sigma $.
\item[$(iv)$] $\log(k\sigma)=k\log \sigma $, for any $k\in Aut(J)$.\\
\end{enumerate}
\end{proposition}

Let $\widetilde{\sigma}=(\sigma,\theta)$ and $\widetilde{\tau}=(\tau,\phi)$ are two elements of $\widetilde{S}$. We say that $\widetilde{\sigma}$ and $\widetilde{\tau}$ are  transversal \index{Transversal} if the projections are transversal, $\sigma\top\tau$. Then there exists $u\in U$ such that $u^{-1}(\tau)=-e$ and $u^{-1}(\sigma)\top -e$. We can apply (\ref{log}) to  $u^{-1}(\sigma)$ and define the {\it Souriau index} of the pair $(\widetilde{\sigma},\widetilde{\tau})$ to be
\begin{equation}\label{Souriau-index}\index{Index!Souriau}
m(\widetilde{\sigma},\widetilde{\tau})=\frac{1}{\pi}\bigl[\frac{1}{i}\tr(\log u^{-1}(\sigma))-r(\theta-\phi)\bigr].
\end{equation}

\begin{theorem}[\cite{Clerc-Koufany}]
The Souriau index is $\mathbb{Z}-$valued continuous function on $S^2_\top$ and is invariant under the action of the covering $\Gamma$.
\end{theorem}
We also prove the following essential cohomological property
\begin{theorem}[\cite{Clerc-Koufany}]\label{Leray-Formula}
If $\widetilde{\sigma}_1, \widetilde{\sigma}_2, \widetilde{\sigma}_3\in \widetilde{S}$ have pairwise transverse projections $\sigma_1,\sigma_2,\sigma_3$, then
$$m(\widetilde{\sigma}_1,\widetilde{\sigma}_2)+m(\widetilde{\sigma}_2,\widetilde{\sigma}_3)+m(\widetilde{\sigma}_3,\widetilde{\sigma}_1)$$
is an integer and coincides with the triple Maslov index of the triplet $(\sigma_1,\sigma_2,\sigma_3)$,
\begin{equation*}
\imath(\sigma_1,\sigma_2,\sigma_3)=m(\widetilde{\sigma}_1,\widetilde{\sigma}_2)+m(\widetilde{\sigma}_2,\widetilde{\sigma}_3)+m(\widetilde{\sigma}_3,\widetilde{\sigma}_1)
\end{equation*}
\end{theorem}
Following an idea of \cite{deGosson}, we extend the definition of the Souriau index to $S\times S$ and prove Theorem \ref{Leray-Formula} without the transversality condition.\\

Now, fix a Jordan frame $(c_j)_{1\leq j\leq r}$, then we have
\begin{theorem}[\cite{Clerc-Koufany}]
Fix $$\widetilde{\sigma}_1=(\sum_{j=1}^\ell e^{i\theta_j}c_j+\sum_{j=\ell+1}^r e^{i\theta_j}
c_j,\theta),\;\;
\widetilde{\sigma}_2=(\sum_{j=1}^\ell e^{i\theta_j}c_j+\sum_{j=\ell+1}^r e^{i\varphi_j}
c_j,\varphi)\in \widetilde{S}$$ such that  $\theta_j-\varphi_j\notin2\pi\mathbb{Z}$ for
$\ell+1\leq j\leq r$. Then
\begin{equation}
  \label{eq:primitive}
 m(\widetilde{\sigma}_1,\widetilde{\sigma}_2)=\frac{1}{\pi}
\Bigl[\sum_{j=\ell+1}^r\{\theta_j-\varphi_j+\pi\}-r(\theta-\varphi)\Bigr]. 
\end{equation}
\end{theorem}
In particular, if $\widetilde{\sigma}_1=(-e,-\pi)$ and
$\widetilde{\sigma}_2=(-\sum_{j=1}^\ell c_j+\sum_{j=\ell+1}^re^{i\varphi_j}c_j,\varphi)$,
where
\begin{itemize}
\item[$(i)$] $-\pi<\varphi_j<\pi$, $\forall j$, $\ell+1\leq j\leq r$, et
\item[$(ii)$] $r\varphi=-\ell\pi+\sum_{j=\ell+1}^r\varphi_j+2k\pi$, with
  $k\in\mathbb{Z}$,
\end{itemize}
then
\begin{equation}
  \label{eq:m_en_coord_simple}
  m(\widetilde{\sigma}_1,\widetilde{\sigma}_2)=2k+r-\ell=2k+r-\mu(\sigma_1,\sigma_2),\\
\end{equation}

where $\mu(\sigma_1,\sigma_2)$ is the transversality index
\index{Index!transversality} of the pair $(\sigma_1,\sigma_2)$.

\newpage


 \section{The Arnold-Leray index}\label{section11}\index{Index!Arnold-Leray}
{\it In this section we generalize the notion of Maslov cycles, prove their
stratification and their causal orientation. We use these geometric
properties to construct  another primitive to the triple Maslov index. }

 \subsection{The Maslov cycles}\index{Maslov cycle}

 Given $\sigma_0\in S$ we consider the following sub-manifold of point of $S$ which are not transversal to $\sigma_0$,
\begin{equation*}
\Sigma(\sigma_0)=S\setminus S_\top(\sigma_0)=\{\sigma\in S \tq \det(\sigma-\sigma_0)=0\}.
\end{equation*}
This is the {\it Maslov cycle} attached to $\sigma_0$.
\begin{theorem}[\cite{Clerc-Koufany}]
The Maslov cycle $\Sigma(\sigma_0)$ is a stratified sub-manifold of $S$ of codimension $1$ and its singularity is of codimension $\geq 3$ in $S$.
\end{theorem}
In fact the stratums are the sets
$$S_k(\sigma_0)=\{\sigma\in S \; ;\; \mu(\sigma,\sigma_0)=k\},\; 0\leq k\leq r.$$
We prove that $S_k(\sigma_0)$ is a sub-manifold of $S$ of codimension $k+\frac{k(k-1)}{2}d$, and 
\begin{equation*}
\Sigma(\sigma_0)=\bigsqcup_{1\leq k\leq r}S_k(\sigma_0)=\overline{S_1(\sigma_0)}.
\end{equation*}
The singularity of the Maslov cycle $\Sigma(\sigma_0)$ is $\bigsqcup_{2\leq
  k\leq r}S_k(\sigma_0)=\overline{S_2(\sigma_0)}$, and the set of regular
points is $S_1(\sigma_0)$.\\

Now we wish  to prove that the Maslov cycle is in addition oriented. For this purpose we use one more time the causal structure of the Shilov boundary $S$ and prove the following fundamental fact :
\begin{proposition}[\cite{Clerc-Koufany}]
Let $\sigma\in S_1(\sigma_0)$ be a regular point of $\Sigma(\sigma_0)$. Set $H_{\sigma_0}(\sigma)=T_\sigma(\Sigma(\sigma_0))$. Then the tangent vectors of all causal curves starting from $\sigma$ are all contained in the same half space of $T_\sigma(S)$ limited by $H_{\sigma_0}(\sigma)$. 
\end{proposition}
Hence $H_{\sigma_0}(\sigma)$ has two sides : $+$ side and $-$ side.  The $+$ side is the one which contains all mentioned tangent vectors. We will denote it    by $H_{\sigma_0}^+(\sigma)$. The family 
$$(H_{\sigma_0}^+(\sigma))_{\sigma\in S_1(\sigma_0)}$$
is called the {\it canonical transverse orientation} \index{Canonical transverse orientation} of the Maslov cycle $\Sigma(\sigma_0)$.  We now claim that this orientation is compatible with the action of the group $G$ : let $g\in G$, then the transverse orientation of  the Maslov cycle $\Sigma(g\sigma_0)$ is given by the family
$$Dg(\sigma)\left[H^+_{\sigma_0}(\sigma)\right]=H^+_{g\sigma_0}(g\sigma),\;\; \sigma\in S_1(\sigma_0),$$
since the causal structure of $S$ is $G-$invariant.

\subsection{The Arnold-Leray index}
 We now wish to use the topological properties of Maslov cycles to construct a homotopy invariant.\\
Let $\sigma_0\in S$ and set $\Sigma_0=\Sigma(\sigma_0)$. A {\it proper
  path} (relatively to $\Sigma_0$)  \index{Proper path} in $S$ is a smooth
path $\gamma : [0,1]\to S$ such that $\gamma(0)\not\in\Sigma_0$,
$\gamma(1)\not\in\Sigma_0$ and intersects $S_1(\sigma_0)$ transversally in
a finite number of crossings, say in $t_1, t_2, \ldots, t_k$.\\

We define the {\it Arnold number} \index{Arnold number}  $\nu_A(\gamma)$ of the proper path $\gamma$ to be the number of intersections of $\gamma$ with $S_1(\sigma_0)$, each counted with sign $\pm$ according to whether the crossing is in the positive or negative direction. Or in the same thing,
\begin{equation*}
\nu_A(\gamma)=\epsilon_1+\ldots+\epsilon_k,\;\; \textrm{with}\;\;
\epsilon_j=
\begin{cases}
+1 & \text{if $\dot{\gamma}(t_j)\in H^+_{\sigma_0}(\gamma(t_j))$}\\
-1 & \text{if $\dot{\gamma}(t_j)\in H^-_{\sigma_0}(\gamma(t_j))$}
\end{cases}
\end{equation*}
It is easy to show that  Arnold number satisfies the following properties (see \cite{Clerc-Koufany}) :
\begin{itemize}
\item Every homotopy class of a given path contains a proper path with the same endpoints.
\item Two homotopic proper paths with the same endpoints have the same Arnold number.
\end{itemize}
This allows us to define the Arnold number for any path $\gamma$ to be the Arnold number of a proper path $\gamma_{\text{proper}}$ homotopic to $\gamma$ with the same endpoints,
$$\nu_A(\gamma)=\nu_A(\gamma_{\text{proper}}).$$
We will now define an index of a pair of point of the universal covering $\widetilde{S}$ by using the invariance by homotopy of the Arnold number.\\
Let $\widetilde \sigma_0, \widetilde \tau_0\in \widetilde S$, and $\sigma_0, \tau_0\in S$ their projections. Let $\sigma(t)$, $0\leq t\leq 1$ be a causal curve\index{Causal curve} such that $\sigma(0)=\sigma_0$. Let $\tau(t)$, $0\leq t\leq 1$ be an anti-causal curve such that $\tau(0)=\tau_0$. Let $\widetilde{\sigma}(t)$ be the lift of $\sigma(t)$  with origin  $\widetilde{\sigma}_0$, and $\widetilde{\tau}(t)$ be the lift of $\tau(t)$,  with origin  $\widetilde{\tau}_0$.\\

 Then we claim, see \cite{Clerc-Koufany}, that there exists $\epsilon>0$ such that for all $t$, $0<t<\epsilon$ the points $\sigma(t)$ and $\tau(t)$ are outside of Maslov cycle $\Sigma(\sigma_0)$ attached with $\sigma_0$. Fix a such $t$ and let $\gamma_t(s)$, $0\leq s\leq 1$ be a proper path (relatively to $\Sigma(\sigma_0)$) with origin $\sigma(t)$ and  end $\tau(t)$ such that its lift is of origin $\widetilde{\sigma}(t)$ and  end $\widetilde{\tau}(t)$.\\

 We define the {\it Arnold index} \index{Index!Arnold} $\nu(\widetilde{\sigma}_0,\widetilde{\tau}_0)$ of the pair $(\widetilde{\sigma}_0,\widetilde{\tau}_0)$ to be the Arnold number of the path $\gamma_t$,
$$\nu(\widetilde{\sigma}_0,\widetilde{\tau}_0)=\nu_A(\gamma_t).$$


Clearly, for fixed $t$, this index does not depend on the choice of the path $\gamma_t$, since any other path having the same properties is homotopic to $\gamma_t$ and thus has the same Arnold  number.  Moreover we prove that $\nu(\widetilde{\sigma}_0,\widetilde{\tau}_0)$ does not depend on the parameter $t$, and on which of the (causal or anti-causal) curves $\sigma(t)$ and $\tau(t)$ we use.\\

The construction of the Arnold index uses only the invariant concepts by causal transformations (causal curves, Maslov cycles, transversality), and thus the Arnold index is invariant under the action of the group $\Gamma$.\\

We will now calculate the index of Arnold `` in coordinates '', as we did for the index of Souriau.  One fixes for that a Jordan frame $(c_j)_{ 1\leq j\leq r}$ of $V$. Thanks to the results concerning the orbits of the action of $G$ in  $S\times S$ (see section \ref{section_causal}), the following proposition covers the general case.

\begin{proposition}[\cite{Clerc-Koufany}]
 Let  
$$\widetilde{\sigma}_0=\widetilde{-e} =(-e,-\pi)\;\;\text{and}\;\;
\widetilde{\tau}_0=\big(-\sum_{j=1}^\ell c_j+\sum_{j=\ell+1}^r
e^{i\varphi_j}c_j,\ \varphi\big)\in \widetilde{S}.$$
One notes $\sigma_0=-e$ and $\tau_0$ their corresponding projections on $S$. 
Suppose
\begin{enumerate}
\item[$(i)$] $-\pi<\varphi_j<\pi$, $\forall j$, $\ell+1\leq j\leq r$;
\item[$(ii)$] $r\varphi= -\ell\pi+\sum_{j=\ell+1}^r\varphi_j+2k\pi$,
  with $k\in \mathbb{Z}.$
\end{enumerate}
Then 
\begin{equation}\label{Indice_Arnold_en_Coor}
\nu(\widetilde{\sigma_0},{\widetilde \tau}_0)
=-\ell+k=k-\mu(\sigma_0,\tau_0).
\end{equation}
\end{proposition}

\begin{corollary}[\cite{Clerc-Koufany}] 
\begin{equation}\label{lien_m_nu}
\nu(\widetilde{\sigma},\widetilde{\tau})=\frac{1}{2}\bigl[m(\widetilde{\sigma},\widetilde{\tau})-\mu(\sigma,\tau)-r\bigr].
\end{equation}

\end{corollary}
A consequence of this corollary is that the right-hand side of the formula (\ref{lien_m_nu}) is an integer. This allows us to introduce the {\it index of inertia} and the\index{Index!of inertia} {\it Arnold-Leray} index.\index{Index!Arnold-Leray}\\

Let $(\sigma_1,\sigma_2,\sigma_3)\in S^3$.  We define  the inertia index of the triplet 
$(\sigma_1,\sigma_2,\sigma_3)$  to be  
\begin{equation*} 
\jmath(\sigma_1,\sigma_2,\sigma_3) = \frac{1}{2}\big(\iota(\sigma_1,\sigma_2,\sigma_3)+\mu(\sigma_1,\sigma_2) -\mu(\sigma_1,\sigma_3)+\mu(\sigma_2,\sigma_3)+r\big).
\end{equation*}
 Let  $\widetilde{\sigma}_1,\widetilde{\sigma}_2 \in
    \widetilde{S}$, and    $\sigma_1,\sigma_2$ the corresponding projections. We define the  Arnold-Leray index to be  
\begin{equation*}
n(\widetilde{\sigma}_1,\widetilde{\sigma}_2) = \nu(\widetilde{\sigma}_1,\widetilde{\sigma}_2)+\mu(\sigma_1,\sigma_2)+r.
\end{equation*}
 One finishes this section by announcing this theorem
\begin{theorem}[\cite{Clerc-Koufany}] 
The index of inertia satisfies  the following :
\begin{itemize}
\item[$(i)$] $\jmath$ is $\mathbb{Z}-$valued function.
\item[$(ii)$] $\jmath$ is a 2-cocycle\footnote{It will be observed that the index  of inertia $\jmath$ does not satisfy the skew-symmetric property that  has the triple Maslov index $\imath$.}
\begin{equation*}
\jmath(\sigma_1,\sigma_2,\sigma_3)-\jmath(\sigma_1,\sigma_2,\sigma_4)+\jmath(\sigma_1,\sigma_3,\sigma_4)-\jmath(\sigma_2,\sigma_3,\sigma_4)=0
\end{equation*}
fro all  $\sigma_1,\sigma_2,\sigma_3\in S$.
\item[$(iii)$] The Arnold-Leray index is a  primitive of the index of inertia, i.e.
\begin{equation*}
\jmath(\sigma_1,\sigma_2,\sigma_3) =n(\widetilde{\sigma}_1,\widetilde{\sigma}_2)-n(\widetilde{\sigma}_1,\widetilde{\sigma}_3)+n(\widetilde{\sigma}_2,\widetilde{\sigma}_3)
\end{equation*}
for all  $\widetilde{\sigma}_1,\widetilde{\sigma}_2,\widetilde{\sigma}_3\in \widetilde{S}$, with the corresponding projections  $\sigma_1,\sigma_2,\sigma_3$.
\end{itemize}
\end{theorem}

\newpage

\section[The Poincar{\'e} rotation number]{The Poincar{\'e} rotation number of the conformal group}\label{extra_section}
{\sl In This section we use the Souriau index to generalize the notion of
  Poincar{\'e} rotation number}\\

A group $G$ is {\it uniformly perfect} \index{Group!uniformly perfect}, if
there exists an integer $k$ such that, every element $g\in G$ is a product
of $k$ commutators at the maximum. A such group has the following
property :\\

{\it 
Let $G$ be a uniformly perfect group and let 
\begin{equation*}
 \begin{CD}
  0 @> >> \mathbb{Z} @>\iota>> \Gamma @>\pi >> G @> >> 1
 \end{CD} 
\end{equation*}
be a central extension of $G$. Let $T=\iota(1)$. Then there exists at the
maximum one map $\Phi : \Gamma \to \mathbb{R}$ such that
\begin{itemize}
\item[$(1)$] $\Phi(\gamma T) = \Phi(\gamma) +1$, $\forall \gamma \in \Gamma$
\item[$(2)$] $ \Phi(\gamma_1 \gamma_2)-\Phi(\gamma_1)-\Phi(\gamma_2)$
is bounded on $\Gamma\times \Gamma$
\item[$(3)$] $\Phi(\gamma^n)=n\Phi(\gamma)$, $\forall \gamma \in
\Gamma$, $\forall n\in \mathbb{Z}$.\\
\end{itemize}
}
A map $\Phi$ satisfying (2) is called a {\it
  quasi-morphism}\index{quasi-morphism}. A map $\Phi$  satisfying (2) and
(3) is called a {\it homogeneous quasi-morphism}.\\

If $\Phi$ exists, then the function 
$$c : G\times G\to\mathbb{R},\;\;
c(g_1,g_2)=\Phi(\gamma_1\gamma_2)-\Phi(\gamma_1)-\Phi(\gamma_2)$$ (where
$\gamma_i$ is the lift of $g_i$) is well defined and is a 2-cocycle, {\it i.e.\/}
$$c(g_1,g_2)+c(g_1g_2,g_3)=c(g_1,g_2g_3)+c(g_2,g_3).$$
Let us consider the example $G=Homeo^+(S^1)$, the group of all homeomorphisms
of the circle preserving the orientation, where $S^1$ is the oriented unite
circle. This group is uniformly perfect
and we consider the central extension
\begin{equation*}
 \begin{CD}
  0 @> >> \mathbb{Z} @>\iota>> \widetilde{Homeo^+(S^1)} @>\pi >> Homeo^+(S^1) @> >> 1
 \end{CD} 
\end{equation*}
where $\widetilde{Homeo^+(S^1)}$ is the universal covering of
$Homeo^+(S^1)$.
We will exhibit a function satisfying (1), (2) and (3).\\
 Let us introduce
the cyclic order. Let $p,q,r\in S^1$, the {\it cyclic order}\index{cyclic
  order} of $p,q,r$ is defined by
$$\text{ord}(p,q,r)=\begin{cases} 
0 & \text{if 2 points coincide}\\ 
1 & \text{if\;} q\in )p,r(\\
-1& \text{if\;} q\in )r,p(.
 \end{cases}$$
Let
$$\Phi_{\text{ord}} : \widetilde{Homeo^+(S^1)} \to \mathbb{Z}, \;
\Phi_{\text{ord}}(\widetilde{f})=2E(\widetilde{f}(0))$$
where $$E(x)=\begin{cases} 
x &\text{if\;} x\in\mathbb{Z},\\
[x]+\frac{1}{2} &\text{if\;} x\notin\mathbb{Z}.
\end{cases}$$
Then $\Phi_{\text{ord}}$ is (non homogeneous) quasi-morphism. Indeed,
$$\Phi_{\text{ord}}(\widetilde{f}\circ\widetilde{g})-\Phi_{\text{ord}}(\widetilde{f})-\Phi_{\text{ord}}(\widetilde{g})=\text{ord}(1,f(1),f\circ
g(1))$$
which is bounded. Moreover, the following limit exists
$$\lim_{n\to
  +\infty}\frac{1}{n}\Phi_{\text{ord}}(\widetilde{f}^n)=2\lim_{n\to
  +\infty}\frac{E(\widetilde{f}^n(0))}{n}=2\tau(\widetilde{f}).$$
The function $\tau : \widetilde{Homeo^+(S^1)} \to\mathbb{R}$ is the
so-called {\it Poincar{\'e} translation number\/} \index{Poincar{\'e} translation
  number}. This function satisfies (1), (2) and (3). Passing to the
quotient, we get  a function $\rho : Homeo^+(S^1) \to
\mathbb{R}/\mathbb{Z}\simeq S^1$, which is the so-called {\it Poincar{\'e}
  rotation number\/} \index{Poincar{\'e} rotation number}.\\

We return to the general case of Hermitian symmetric spaces of tube type
$\mathcal{D}=G/K$. We will use the triple Maslov index and the Souriau
index to generalize the
notion of the Poincar{\'e} rotation number.  We prove  first the following 

\begin{proposition}[\cite{Clerc-Koufany}] Let $G=KAN$ be the Cartan decomposition of
$G$, then
\begin{itemize}
\item[$(i)$] Every element of $N$ is a commutator.
\item[$(ii)$] Every element of $A$ is a product of $r$ commutators at the maximum.
\end{itemize}
\end{proposition}

We also need the following lemma
\begin{lemma}[\cite{Clerc-Koufany}]  Consider the following central extension of $G$
\begin{equation*}
 \begin{CD}
  0@> >>\mathbb{Z} @> \imath >> \Gamma @>p>> G @> >> 1.
 \end{CD} 
\end{equation*}
Then there exists at the maximum one map $\Phi : \Gamma \to\mathbb{R}$ such
that
\begin{itemize}
\item[$(0)$] $\Phi$ is continuous 
\item[$(1)$] $\Phi(\gamma T) = \Phi(\gamma) +1$, $\forall \gamma \in \Gamma$
\item[$(2)$] $ \Phi(\gamma_1 \gamma_2)-\Phi(\gamma_1)-\Phi(\gamma_2)$
is bounded on $\Gamma\times \Gamma$
\item[$(3)$] $\Phi(\gamma^n)=n\Phi(\gamma)$, $\forall \gamma \in
\Gamma$, $\forall n\in \mathbb{Z}$.
\end{itemize}
\end{lemma}

We begin now to construct the Poincar{\'e} rotation number on $G$. Let
$\widetilde{o}$ be a base point of $\widetilde{S}$, the universal covering of
the Shilov boundary. Then the function 
$$\c : \Gamma \to \mathbb{Z}  \;  : \;
\c(\gamma)=m(\gamma\cdot\widetilde{o},\widetilde{o})$$ where $m$ is the
Souriau index, is a (non homogeneous) quasi-morphism. Indeed, 
$$\begin{array}{ll}
\c(\gamma_1\gamma_2)-\c(\gamma_1)-\c(\gamma_2)&=m(\gamma_1\gamma_2\cdot\widetilde{o},\widetilde{o})+m(\widetilde{o},
\gamma_1\cdot\widetilde{o})+ m(\gamma_1\cdot
\widetilde{o},\gamma_1\gamma_2\cdot\widetilde{o})\\
&=\imath(o,g_1\cdot o,g_1g_2\cdot o),
\end{array}
 $$
which is bounded by the rank $r$,
where  $\gamma_j$ is the lift of $g_j$, $j=1, 2, 3$.\\
Hence, for $\gamma\in\Gamma$, the sequence $\c_k=\c(\gamma^k)$ satisfies
$$|\c_{k+\ell}-\c_k-\c_\ell|\leq r.$$
Thus, the following limit exists
$$\lim_{k\to+\infty}\frac{1}{k}\c(\gamma^k):=\tau(\gamma).$$

\begin{theorem}[\cite{Clerc-Koufany}]
The function $-\frac{1}{2}\tau$ is a continuous homogeneous quasi-morphism of $\Gamma$
and it is independent of the choice of $\widetilde{o}$. 
\end{theorem}
The function $\tau$ is the {\it generalized Poincar{\'e} translation number}
\index{Generalized!Poincar{\'e} translation number}. Passing to the quotient,
the function
$$\rho(g) = -\frac{1}{2}\tau(\gamma) \mod(\mathbb{Z})$$
where $\gamma$ is the lift of $g$, is the {\it generalized Poincar{\'e}
  rotation number}\index{Generalized!Poincar{\'e} rotation number} of
$G$. Finally, we prove the following

\begin{proposition}[\cite{Clerc-Koufany}]The function $\rho$ satisfies ;
\begin{enumerate} 
\item $\rho$ is invariant by conjugaison.
\item If $g\in G$ fixes a point in $S$, then $\rho(g)=0$.
\item If $u\in U$, then $e^{2i\pi\rho(u)}=\chi(u)$
\end{enumerate}
\end{proposition}

\part{Non-commutative Hardy spaces }

\section{Hardy spaces on Lie semi-groups}\label{section12}\index{Hardy space}

{\it In this section we recall the theory of Hardy spaces on Lie semigroups
  due to Olshanski\u{\i} }\\

Let $\frak g$ be a simple Lie algebra over the reals $\Bbb R$, and $\frak
g=\frak k\oplus\frak p$ a Cartan decomposition of $\frak g$. Let $\frak
t\subset\frak k$ be a Cartan subalgebra of $\frak k$. We shall suppose
that $\frak k$ has a non--zero center $\frak z$ ; then $\frak z$ is one
dimensional and $\frak t$ is also a Cartan subalgebra of $\frak g$. 

Let $G_{\Bbb C}$ be the simply connected complex Lie group corresponding
to $\frak g_{\Bbb C}:=\frak g+i\frak g$, and let $G$, $K$ and $T$ be the connected subgroups in $G_{\Bbb C}$
corresponding to $\frak g$, $\frak k$ and $\frak t$ respectively. By the Kostant--Paneitz--Vinberg Theorem 
\cite{Vin}, there are non--trivial regular cones $C$ in
$i\frak g$ which are $\mathrm{Ad}(G)-$invariant, where regular means, convex,
closed, pointed ($C\cap-C=\{0\}$) \index{Cone!pointed} and generating ($C-C=i\frak g$).\index{Cone!generating} Let
$\mathrm{Cone}(i\frak g)$ be the set of all regular $\mathrm{Ad}(G)-$invariant cones in
$i\frak g$. 

For such a cone $C$ in $\mathrm{Cone}(i\frak g)$, Ol'shanski\u{\i} associates a
semigroup $\Gamma(C):=G\exp(C)$ in $G_{\Bbb C}$, and for this semigroup
he associates a ``non-commutative'' Hardy space $H^2(\Gamma(C))$ which is
the set of holomorphic functions $f$ on the complex manifold
$\Gamma(C)^\circ=G\exp(C^\circ)$, the interior of $\Gamma(C)$, such that
$$\sup_{\gamma\in \Gamma(C)^\circ}\int_G |f(g\gamma)|^2\,dg<\infty.$$
For any $\gamma\in \Gamma(C)^\circ$ the linear functional $f\longmapsto
f(\gamma)$ is continuous on $H^2(\Gamma(C))$. Therefore by the Riesz
representation theorem, there exists a vector $K_\gamma\in
H^2(\Gamma(C))$ such that $(f , K_\gamma)=f(\gamma)$. The reproducing
kernel $K$ which is called the {\it Cauchy--Szeg{\"o} kernel} \index{Cauchy--Szeg{\"o} kernel} is defined
by 
$$K(\gamma_1 , \gamma_2)=K_{\gamma_2}(\gamma_1).$$
It is Hermitian, holomorphic in $\gamma_1$ and anti--holomorphic in
$\gamma_2$.

Let
$\Delta=\Delta(\frak g_{\Bbb C}, \frak t_{\Bbb C})$ be the set of roots \index{Roots} of
$\frak g_{\Bbb C}$ relative to $\frak t_{\Bbb C}$. 
Let $\Delta^+\subset\Delta$ be the set
of positive roots\index{Roots!positive}
relative to some order (namely the one where the
center of $\frak k$ comes first), $\Delta^+_{\frak k}$ and
$\Delta^+_{\frak p}$ the set of positive compact \index{Roots!compact} and non-compact
roots \index{Roots!non-compact}, respectively.
Put $\frak t_{\Bbb R}:=i\frak
t\subset\frak t_{\Bbb C}$. 
We identify $\frak t_{\Bbb R}$ with its own dual via the Cartan--Killing form. Then we can consider $\Delta\subset\frak t_{\Bbb
R}$. 
Let ${\mathcal P}\subset \frak t_{\Bbb R}^*\simeq \frak t_{\Bbb R}$ be the set of weights relative to $T$ and let ${\mathcal R}$ be the set of all highest weights relative to
$\Delta^+_{\frak k}$,
\begin{equation*}
{\mathcal R}=\{\lambda\in {\mathcal P} \tq (\forall \alpha\in \Delta^+_{\frak k}) \,\, \langle \lambda , \alpha
\rangle \geq0 \}.
\end{equation*}
Let $\rho$ be the half sum of all positive
roots. Then by Harish--Chandra (\cite{HC1},\cite{HC2},\cite{HC3}) the holomorphic discrete series representations for the group $G$
are  those
irreducible unitary representations of $G$ that are square--integrable
with a highest weight $\lambda$ belonging to 
\begin{equation*}
{\mathcal R}'=\{\lambda\in {\mathcal R} \tq (\forall \beta\in \Delta^+_{\frak p})\,\, \langle \lambda+\rho , \beta \rangle < 0 \}.
\end{equation*}


We will say that $\lambda\in {\mathcal R}$ satisfies the {\it Harish--Chandra condition} if 
\begin{equation*}
\langle \lambda+\rho,\beta\rangle <0, \quad \forall \beta\in \Delta^+_{\frak p}.
\end{equation*}

By Vinberg \cite{Vin}, there exists in $\mathrm{Cone}(i\frak g)$ a unique (up to multiplication by $-1$) maximal cone \index{Cone!maximal} $C_{\rm max}$, such that
$$C_{\rm max}\cap\frak t_{\Bbb R}=c_{\rm max}:=\{X\in \frak t_{\Bbb R} \tq (\forall\alpha\in\Delta^+_{\frak p})\,\, \langle X,\alpha\rangle\geq0\},$$
and a unique minimal cone $C_{\rm min}=C^*_{\rm max}$, such that $C_{\rm min}\cap\frak t_{\Bbb R}=c_{\rm min}$ is the convex cone spanned by all $\alpha$ in $\Delta^+_{\frak p}$. 
\bigskip
A unitary representation $\pi$ of $G$ in a Hilbert space ${\mathcal H}$ is said to be
$C-${\it dissipative} \index{Dissipative representation} if for all $X\in C$ and all
$\xi\in{\mathcal H}^{\infty}$, the space of ${\mathcal C}^\infty$ vectors in ${\mathcal
H}$,
$$(\pi(X)\xi|\xi)\leq0.$$ 
We can now state the Theorem B of Ol'shanski\u{\i} \cite{Ol2} on the non--commutative
Hardy spaces 

\begin{theorem}[\cite{Ol2}]
The Hardy space $H^2(\Gamma(C))$ is a non--trivial Hilbert space for any $C\in
\mathrm{Cone}(i\frak g)$.\par\noindent
The representation of $G$ in $H^2(\Gamma(C))$ can be decomposed
into a direct sum of irreducible unitary representations of $G$. The
components of this decomposition are precisely all the holomorphic
discrete series representations of $G$ which are $C-$dissipative.
\end{theorem}

The group $G\times G$ acts on $H^2(\Gamma(C))$ via left and right
regular representations. Therefore
\begin{equation}\label{1.4}
H^2(\Gamma(C))=\bigoplus_{\lambda\in(C^*\cap\frak t_{\Bbb R})\cap{\mathcal R}'}
\pi_\lambda\otimes \pi_\lambda^*,
\end{equation}
where $\pi_\lambda$ is the {\it contraction representation} of
$\Gamma(C)$ corresponding to a unitary highest weight representation of
$G$ with highest weight $\lambda$. Moreover, the corresponding function
of the Cauchy--Szeg{\"o} kernel $K$ of $H^2(\Gamma(C))$ can be written on
$\Gamma(C)^\circ$ as follows
\begin{equation}\label{1.5}
K(\gamma):=K(\gamma,e)=\sum_{\lambda\in (C^*\cap\frak t_{\Bbb R})\cap{\mathcal R}'} d_\lambda
\tr(\pi_\lambda(\gamma)),
\end{equation}

where $d_\lambda$ denotes the {\it formal dimension} of the
representation $\pi_\lambda$. The series for $K$ converges uniformly on
compact subsets in $\Gamma(C)^\circ$. 

\begin{remark}
Whenever $C$ is the minimal cone, $C_{min}$, the decompositions (\ref{1.4}) and
(\ref{1.5}) are over all the holomorphic discrete series, i.e. over $\lambda\in
{\mathcal R}'$.\\
One of the most important problems in this areas is to give an explicit formula for the function $K(\gamma)$.
\end{remark}

\newpage

\section{The contraction semigroup}\label{section13}\index{Semi-group!of contractions}
{\it In this section we restrict our self to 
  $G=Sp(r,\mathbb{R})$, $G=SO^(2\ell)$ and $G=U(p,q)$. We prove that in
  this case the Lie semigroup is a semigroup of contractions. We also prove
that the image of this semigroup under a new Cayley transform
is the tube domain (modulo some singular points).}\\

From now on we assume that $G$ is one of the classical groups $U(p,q)$,
$Sp(r,\Bbb R)$ or $SO^*(2l)$. Let $\sigma$ be an involution in $\frak
g_{\Bbb C}$ such that 
$$\frak g=\{X\in \frak g_{\Bbb C} \tq
\sigma(X)=-X\}.$$
 Then 
$$\sigma(X)=JX^*J,$$
where $X^*$ is the adjoint matrix and
$$\begin{array}{lll}
\mathrm{for}& \frak g=\frak{u}(p,q),&
               J=\left(\begin{smallmatrix}-I_p&0\\ 0&I_q \end{smallmatrix}\right)\\
\mathrm{for}& \frak g=\frak{sp}(r, \Bbb R),&
               J=\left(\begin{smallmatrix}-I_r&0\\ 0&I_r\end{smallmatrix}\right)\\
\mathrm{for}& \frak g=\frak{o}^*(2l),&
               J=\left(\begin{smallmatrix}-I_l&0\\ 0&I_l\end{smallmatrix}\right).\\
\end{array}$$
\begin{remark}\label{remark2.1}
$U(p,q)$ is not a simple Hermitian Lie groups. Since
$$U(p,q)\simeq \bigl(U(1)\times SU(p,q)\bigr)/{\Bbb Z_{p+q}},$$
 the holomorphic discrete series representations of $U(p,q)$ are
the holomorphic discrete series representations of the circle times
the Hermitian group $SU(p,q)$ which are
trivial on $(\xi,\xi^{-1}I_{n})$, for $\xi^n=1$ ($n=p+q$). Therefore one can easily generalize the results of
section 1 to the reductive group $U(p,q)$.
\end{remark}

Let $C$ be the regular cone in
$i\frak g$ defined by
\begin{equation*}
C:=\{X\in i\frak g \tq JX\leq0\},
\end{equation*}
and let $\Gamma(C):=G\exp(C)$ be the corresponding Ol'shanski\u{\i}
semigroup. An element $\gamma$ of $G_{\Bbb C}$ is said to be a
{\it $J-$contraction} (resp. a {\it strict $J-$contraction}) if
$J-\gamma^*J\gamma\geq0$ (resp. $J-\gamma^*J\gamma\gg0$).


\begin{proposition}[{\cite{Koufany-Orsted3}}]
The semigroup $\Gamma(C)$ is 
the $J-${\it contractions} semigroup,
$$\Gamma(C)=\{\gamma\in G_{\Bbb C} \tq J-\gamma^*J\gamma \geq0\},$$
and $\Gamma(C)^\circ$ is the semigroup of strict $J-$contractions,
$$\Gamma(C)^\circ=\{\gamma\in G_{\Bbb C} \tq J-\gamma^*J\gamma\gg 0\}.$$
\end{proposition}


Now, let $V$ be one of the {\it Jordan algebras} \index{Jordan algebra} $Herm(n, \Bbb C)$, $Sym(2r, \Bbb R)$
or $Herm(l, \Bbb H)$ and let $\Omega$ be the corresponding {\it symmetric
cone}.  Then $\Omega=V^+$ is the set of positive definite matrices in $V$. The tube domain $T_\Omega:=V+i\Omega$ is a
Hermitian symmetric space isomorphic to $G^\flat/K^\flat$, where
$G^\flat$ is $SU(n,n)$, $Sp(2r, \Bbb R)$ or $SO^*(4l)$ respectively and
$K^\flat$ the corresponding maximal compact subgroup, i.e. $S\bigl(U(n)\times U(n)\bigr)$, $U(2r)$ or $U(2l)$ respectively.

\bigskip

Let $\C$ be the {\it Cayley transform} \index{Cayley transform}defined by
\begin{equation*}
\C(Z):=(Z-iJ)(Z+iJ)^{-1}
\end{equation*}
whenever the matrix $(Z+iJ)$ is invertible.


\begin{proposition}[{\cite{Koufany-Orsted3}}]\label{prop2.3}
The Cayley transform $\C$ is a biholomorphic bijection from an open subset
of the tube domain $T_\Omega$ onto the complex manifold
$\Gamma(C)^\circ$. More precisely, if\/ $\Sigma$ denotes the hypersurface $\Sigma=\{Z\in T_\Omega \tq \det(Z+iJ)=0\}$, then 
$$\C(T_\Omega\setminus\Sigma)=\Gamma(C)^\circ.\leqno(2.3)$$
\end{proposition}
Here ``$\det$'' denotes the determinant of the Jordan algebra $V_{\Bbb C}$ (see \cite{Faraut-Koranyi}).\\

\newpage

\section{The holomorphic discrete series}\label{section14}\index{Holomorphic discrete
  series}
{\it In this section we recall the holomorphic discrete series
  representations of $G^\flat$ and explain our strategy to compare the Hardy
space of the Lie semigroup and the Hardy space of the tube domain.}\\

Let $N$ and $R$ be the dimension and the rank of the Jordan algebra $V$. For a complex manifold ${\mathcal M}$ we denote by ${\mathcal O}({\mathcal M})$ the space of holomorphic functions on ${\mathcal M}$.

The group $G^\flat$ acts on $T_\Omega$ via
$$g\cdot Z=(AZ+B)(CZ+D)^{-1}, \quad g=\left(\begin{smallmatrix}A&B\\ C&D\end{smallmatrix}\right),$$
and the scalar--valued holomorphic discrete series representations of
$G^\flat$ are
$$\bigl(U_\lambda(g)f\bigr)(Z)=\det(CZ+D)^{-\lambda}f(g^{-1}\cdot Z),\quad g^{-1}=\left(\begin{smallmatrix}A&B\\ C&D\end{smallmatrix}\right)$$
for $\lambda\geq 2{N\over R}$, which all are unitary and irreducible in the Hilbert spaces
$${\mathcal H}_{\lambda}(T_\Omega):=\{f\in {\mathcal O}(T_\Omega) \tq \int_{T_\Omega}
|f(X+iY)|^2\det(Y)^{\lambda-2{N\over R}}dXdY<\infty\}.$$
Moreover the 
reproducing kernel of ${\mathcal H}_\lambda(T_\Omega)$ is given by
$$K^{T_\Omega}_\lambda(Z,W)=\det\Bigl({Z-W^*\over{2i}}\Bigr)^{-\lambda}.$$


The {\it classical Hardy space} $H^2(T_\Omega)$ on $T_\Omega$ is defined as the space of holomorphic functions $f$ on $T_\Omega$ such that
$$\sup_{Y\in \Omega}\int_V|f(X+iY)|^2\,dX<\infty.$$


\begin{proposition}
The Hardy space $H^2(T_\Omega)$ may be thought of as the space ${\mathcal h}_\lambda(T_\Omega)$ for $\lambda={N\over R}$, and the Cauchy--Szeg{\"o} kernel of $T_\Omega$ is given by
\begin{equation*}
K(Z,W)=\det\Bigl({Z-W^*\over 2i}\Bigr)^{-N/R}.
\end{equation*}
\end{proposition}


We list here the groups $G$ and the corresponding group $G^\flat$,
Jordan algebra $V$, its rank $R$, its dimension $N$, and the determinant $\det$
:
 \begin{center}
\begin{table}[h]
\begin{tabular}[t]{cccccc}
\hline
$G$           & $G^\flat$       & $V$              & $N$       & $R$  & $\det$ \\
\hline
\hline
$Sp(r,\Bbb R)$ & $Sp(2r, \Bbb R)$ & $Sym(2r, \Bbb R)$ & $r(2r+1)$ & $2r$ & $\Det$ \\
${SO}^*(2l)$  & ${SO}^*(4l)$    & $Herm(l,\Bbb H)$  & $l(2l-1)$ & $l$  & $\Det^{1/2}$ \\
$U(p,q)$      & $SU(n,n)$       & $Herm(n,\Bbb C)$  & $n^2$     & $n$  & $\Det$ \\
\hline
\hline
\end{tabular}
\end{table}
\end{center}

A crucial point is to compare holomorphic functions on the tube
domain with their pull-backs on the semigroup via the Cayley transform, and vice
versa. In particular, it will be important to know the rate of growth of
the functions near the singularity $\Sigma$ above. Assuming $\gamma = \C(Z)$
we have that
$$Z+iJ = 2(I-\gamma)^{-1}iJ$$
so that to approach the singularity in the $Z$ variable, means that
$\det(I-\gamma)$ tends to infinity in the $\gamma$ variable. Clearly 
this condition is invariant under conjugation with $G$, so we
may reduce the question of the growth near the singularity
to a question on the compact Cartan subspace.
Suppose the holomorphic functions $f$ and $F$ are related by
$$ f(Z) = \det(I-\gamma)^pF(\gamma) $$
so that $F$ is holomorphic on $\Gamma(C)^\circ$ and $f$
therefore holomorphic on $T_\Omega\setminus\Sigma$. Then for $f$
to admit a holomorphic continuation to all of $T_\Omega$ it is
necessary and sufficient that it stays bounded as the
determinant factor tends to infinity, i.e. that $F$ satisfies
a decay condition related to $p$. This is what we shall make
precise in the following. 

\newpage

\section{The case of $G=Sp(r,\mathbb{R})$}\label{section15}
{\it In the section we will give an explicit construction of double
  covering of the Lie semigroup and compare the two Hardy spaces. In the 
  case $G=Sp(r,\mathbb{R})$, the classical Hardy
space is isomorphic to the odd part of the Olshanski\u{\i} Hardy space}\\

We assume that $G=Sp(r,\Bbb R)$. Then the Hardy parameter is ${N\over R}=r+{1\over2}\in \Bbb Z+{1\over2}$ 
and $\det$ coincides with the usual matrix determinant Det. This suggests that the 
operator 
$\C_{N\over R} = \C_{r+{1\over2}}$, 
\begin{equation}\label{Cayley}
f=\C_{r+{1\over2}}(F) \,\,\, : \,\,\,f(Z)=\Det(Z+iJ)^{-(r+{1\over2})}F(\gamma)
\end{equation}
may be an intertwining operator between $H^2(T_\Omega)$ and the odd part of the Hardy space $H^2(\Gamma(C)_2)$ 
on the double covering $\Gamma(C)_2$ of $\Gamma(C)$.\\
\subsection{The explicit construction of the two sheeted covering semigroup}
For the open subset $\Gamma(C)_2\hskip-1pt^\circ$ we have a new and explicit construction:
Let $J^\flat=\left(\begin{smallmatrix}J&0\\ 0&J \end{smallmatrix}\right)$ and let
$G^\flat:=Sp(2r,\Bbb R)$ be the group of all matrices in $Sp(2r,\Bbb C)$ satisfying
$$g^*J^\flat g=J^\flat.$$
We imbed $G$ in a natural way in $G^\flat$ as follows :
$$g\longmapsto \begin{pmatrix}g&0\cr 0&I_{2r}\end{pmatrix}.$$
We also view the Cayley transform $\C$ as the element of $G^\flat_{\Bbb C}$ given by the matrix
$$\C={1\over\sqrt{2}}\begin{pmatrix}\hfill I_{2r}&-\hfill iJ\cr \hfill
I_{2r}&\hfill iJ\end{pmatrix}.$$
Our precise definition of $G_2$ is to be the set of all pairs  
$(g,\omega(g^{\c},\cdot))$
 with $g\in Sp(r,\Bbb R)$,
$g^{\c}=\C^{-1}\,g\,\C=\left(\begin{smallmatrix}A&B\cr C&D\end{smallmatrix}\right)$,
$\omega(g^{\c},\cdot)^2=\Det(C\,\cdot+D)^{-1}$ and
$Z\longrightarrow \omega(g^{\c},Z)$ is holomorphic on $T_\Omega$.
Note that this is analogous to the definition of the double cover of $SU(1,1)$, 
where we take all pairs
$(g,\sqrt{cz+d})$ with $g=\begin{pmatrix}a&b\cr c&d\end{pmatrix}\in SU(1,1)$ and $\sqrt{cz+d}$ a
 holomorphic
choice of square root of the non-zero function $cz+d$ on the unit disc.
Indeed, it sometimes is convenient to think in terms of such multivalued
functions when doing practical calculations, but of course, the precise
definition is behind this. We also recall
the more informal definition of $G_2$ as follows: \\
Take again $Z\in T_{\Omega}$ and $g\in Sp(r,\Bbb R)$ such that
$\C^{-1}\,g\,\C=\left(\begin{smallmatrix}A&B\cr C&D\end{smallmatrix}\right)$. A determination on
$T_{\Omega}$ of the square root $\Det(CZ+D)^{-{1\over2}}$ is completely
determined by its value on $Z=iI$. For each $g\in Sp(n,\Bbb R)$ we
choose a determination of $\Det(CZ+D)^{-{1\over2}}$. This is a global determination.
We consider here $Z$ as a variable, since the group
(and indeed all contractions) acts on the tube domain, and we consider the function 
\begin{equation*}
Z\longmapsto \omega(g^{\c},Z):=\Det(CZ+D)^{-{1\over2}}
\end{equation*}
from $T_{\Omega}$ into $\Bbb C\setminus\{0\}$, where $g^{\c}=\C^{-1}\,g\,\C$. 
We read $\omega$ as ``a holomorphic choice of square root of the determinant". It follows that $\omega$ may be viewed as a cocycle for $G_2$, and it gives a choice of square root at the
product of two elements as follows:  
$$\omega(g^{\c}_1g^{\c}_2,Z)=\omega(g^{\c}_1,g^{\c}_2\cdot Z)\omega(g^{\c}_2,Z).$$
This equation is to be understood as an equation for the two-valued
function $\omega$; it does not hold for any single-valued function. 
More generally, assuming the determinant to be non-zero, we let
$$\omega_2(\left(\begin{smallmatrix}A&B\cr C&D\end{smallmatrix}\right),Z):=\Det(CZ+D)^{-1}$$
and correspondingly $\omega(\left(\begin{smallmatrix}A&B\cr C&D\end{smallmatrix}\right),Z)$
a choice of one of the two square roots of this, either (as here) global
and holomorphic, or (as below) local, i.e. at the fixed point $Z$.  
Then we may consider our double covering group to be  
$$G_2:=\{\widetilde{g}:=(g,\omega(g^{\c},\cdot))\tq g\in Sp(r,\Bbb R)\},$$
endowed with the group law
\begin{equation}\label{gp_law}
\left(g_1,\omega(g_1^{\c},Z)\right)\left(g_2,\omega(g_2^{\c},Z)\right)=\left(g_1g_2,\omega(g_1^{\c},g_2^{\c}\cdot Z)\omega(g_2^{\c},Z)\right).
\end{equation}
$G_2$ is a two--sheeted covering group of $G$, since we are considering both choices of square root.
$G_2$ is called the {\it metaplectic group}. \index{Group!metaplectic}
 
Now we wish to give another version of the double covering construction. Here $N$ will be the open
semigroup, realized as a subset of the tube domain as in Proposition \ref{prop2.3}. 
For $Z\in T_\Omega\setminus\Sigma$ and for a
choice of a local determination of $\Det(Z+iJ)^{-{1\over2}}$ we note that up to a constant  
$$\Det(Z+iJ)^{-{1\over2}}=\omega(\C,Z).$$
This is again an identity between two-valued functions. 
Hence at each fixed point $Z$ we make a choice between the two possible values
of the square root, so here the notation does not consider $Z$ as a variable.
Note that we may extend our cocycle to the complexified group in the natural way.
Therefore, the complex manifold
$$\Gamma(C)^\circ_2:=\{\widetilde{\gamma}=\left(\gamma,\omega(\C,Z)\right)
 ; \gamma\in \Gamma(C)^\circ,\,\gamma=\C(Z), \, Z\in T_\Omega\setminus\Sigma\},$$
is a two--sheeted covering of the semigroup $\Gamma(C)^\circ$. As before, we consider both choices
 of square root here, and corresponding to the modern point of view, the more
precise definition of $\Gamma(C)^\circ_2$ is the set of $\widetilde{\gamma}=(\gamma,w)\in\Gamma(C)^\circ\times \Bbb C$ such that 
$$
\gamma=\C(Z),\;  Z\in T_\Omega\setminus\Sigma\; \textrm{and} \; w^2=\Det(Z+iJ)^{-1}.$$
In particular, $w$ is just a complex number.

\begin{lemma}[{\cite{Koufany-Orsted3}}]
The group $G_2$ acts on the right on the manifold $\Gamma(C)^\circ_2$
\end{lemma} 

Indeed, letting  $Z'$ satisfy $g^{-1}\gamma = \C\cdot Z'$,  
which implies that $ Z = g^{\c}\cdot Z'$, 
then
$$
\bigl(\gamma,\omega(\C,Z)\bigr)\cdot\left(g,\omega(g^{\c},\cdot)\right)=\bigl(g^{-1}\gamma,\omega(\C,Z')\bigr).$$

To show that $\Gamma(C)^\circ_2$ is a semigroup we consider the following manifold
$${\Gamma(C)^\circ_2}':=
\{\left(\gamma,\omega(\gamma^{\c},\cdot)\right)\tq \gamma\in \Gamma(C)^\circ\}.$$
It is clear that ${\Gamma(C)^\circ_2}'$ is a double covering of $\Gamma(C)^\circ$
 and has a semigroup structure with respect to the law (\ref{gp_law}).

Consider the map $\varphi$ from ${\Gamma(C)^\circ_2}'$ to ${\Gamma(C)^\circ_2}$ defined by
$$\bigl(\gamma,\omega(\gamma^{\c},Z)\bigr)\longmapsto\bigl(\gamma,\omega(\C,Z)\bigr), \quad\hbox{\rm where}\,\,\,\gamma=\C(Z)\in \Gamma(C)^\circ. $$

\begin{lemma}[{\cite{Koufany-Orsted3}}]
$\varphi$ is a homeomorphism from ${\Gamma(C)^\circ_2}'$ onto
${\Gamma(C)^\circ_2}$.
\end{lemma}

\begin{remark}
The semigroup $\Gamma(C)^\circ_2$ is isomorphic to the interior of the {\it metaplectic semigroup} \index{Semi-group!metaplectic} or the {\it Howe oscillator semigroup}. We call it the {\it open metaplectic semigroup}. 
\end{remark}
\subsection{The Hardy space on $\Gamma(C)_2$}

The Hardy space $H^2(\Gamma(C)_2)$ on the {\it metaplectic semigroup} $\Gamma(C)_2$ is the space of holomorphic functions $F\in {\mathcal O}(\Gamma(C)^\circ_2)$ such that
$$\sup_{\widetilde{\gamma}\in \Gamma(C)^\circ_2}\int_{G_2}|F(\widetilde{\gamma}\,\widetilde{g})|^2\,d\,\widetilde{g}<\infty.$$

The compact maximal subgroup $K$ of $G=Sp(r,\Bbb R)$ is isomorphic to $U(r)$ and the maximal split Abelian subalgebra
$$\frak t_{\Bbb R}=\Big\{\left(\begin{smallmatrix}
                        \hfill X&\hfill 0\cr
                        0&-X\end{smallmatrix}
\right)\in {\mathcal M}(r\times r,\Bbb R) \tq 
X=\left(\begin{smallmatrix}x_1 &        &\cr 
                         & \ddots &\cr 
                         &        &x_r\end{smallmatrix}\right)
\Bigr\},$$
can be identified with $\Bbb R^r$.
Let $\epsilon_1, \ldots, \epsilon_r$ be the canonical basis of $\frak
t_{\Bbb R}^*=\frak t_{\Bbb R}=\Bbb R^r$. Then the root system
$\Delta=\Delta(\frak g_{\Bbb C},\frak t_{\Bbb C})$ is of type $C_r$ :
\begin{eqnarray*}
\Delta&=&\{\pm(\epsilon_i\pm\epsilon_j) \,\,\, (1\leq i<j\leq r) \,\,\, , \,\,\, \pm 2\epsilon_i \,\,\, (1\leq i \leq r)\}\\
\Delta^+&=&\{\epsilon_i\pm\epsilon_j \,\,\, (1\leq i<j\leq r) \,\,\, , \,\,\,
2\epsilon_i \,\,\, (1\leq i \leq r)\}\\
\Delta^+_{\frak k}&=&\{\epsilon_i-\epsilon_j \,\,\, (1\leq i<j\leq r)\}\\
\Delta^+_{\frak p}&=&\{\epsilon_i+\epsilon_j \,\,\, (1\leq i\leq j\leq r)\}\\
\rho&=&r\epsilon_1+(r-1)\epsilon_2+\ldots+\epsilon_r\cr
&\simeq&(r, r-1, \ldots , 1).
\end{eqnarray*}
Furthermore ${\mathcal P}$ is the lattice $\Bbb Z^r$, the set of highest weights relative to $\Delta^+_{\frak k}$ is given by
$${\mathcal R}=\{\lambda=(\lambda_1, \ldots, \lambda_r)\in \Bbb Z^r \tq \lambda_1\geq
\ldots\geq\lambda_r\},$$
and $\lambda\in\mathcal R$ satisfies the Harish--Chandra condition if
$$-r>\lambda_1\geq\ldots\geq\lambda_r,$$
which gives the set $\mathcal R'$.

Let $K_2\subset G_2$, resp. $T_2\subset K_2$ be the corresponding covering of $K$ and $T$. Then the corresponding ${\mathcal P}_2$, $\mathcal R_2$ and $\mathcal R'_2$ are given by
\begin{eqnarray*}
{\mathcal P}_2&=&\Bbb Z^r\cup(\Bbb Z^r+{1\over2})={\mathcal P}\cup ({\mathcal P}+{1\over2})={\mathcal P}_{2,\rm even}\cup{\mathcal P}_{2,\rm odd},\\
{\mathcal R}_2&=&\{\lambda\in {\mathcal P}_2 \tq (\forall \alpha\in \Delta^+_{\frak k})\,\, \langle \lambda , \alpha \rangle \geq 0\}\cr
&=&\{\lambda=(\lambda_1, \ldots, \lambda_r)\in\Bbb Z^r\cup(\Bbb
Z^r+{1\over2}) \tq \lambda_1\geq\ldots\geq\lambda_r\},\\
&=&{\mathcal R}_{2,\rm even}\cup {\mathcal R}_{2,\rm odd}\cr
\mathcal R'_2&=&\{\lambda\in{\mathcal R}_2 \tq \langle \lambda+\rho , \beta \rangle<0, \,\,\, \forall \beta\in\Delta^+_{\frak p}\},\\
&=&\{\lambda=(\lambda_1, \ldots, \lambda_r)\in\Bbb Z^r\cup(\Bbb
Z^r+{1\over2}) \tq -r>\lambda_1\geq\ldots\geq\lambda_r\},\\
&=&{\mathcal R}'_{2,\rm even}\cup {\mathcal R}'_{2,\rm odd},
\end{eqnarray*}
where ${1\over2}$ stands for the tuple $({1\over2}, \ldots ,{1\over2})$.
The holomorphic discrete series representations for the metaplectic
group $G_2$ are those irreducible unitary representations
$\pi_\lambda$ of $G_2$ that are square--integrable with a highest
weight $\lambda\in {\mathcal R}'_2={\mathcal R}'_{2,\rm even}\cup {\mathcal R}'_{2,\rm
odd}$. Therefore 
$$H^2(\Gamma(C)_2)=\bigoplus_{
                             \lambda\in(C^*\cap\frak t_{\Bbb R})\cap{\mathcal R}'_2
                            }
                            \pi_{\lambda}\otimes
                            \pi_{\lambda}^*.\leqno(5.4)$$
The cone $C$ is the minimal one in $i\frak g$, so the above
summation is over ${\mathcal R}'_2$ and 
the Hardy space $H^2(\Gamma(C)_2)$ splits into two parts, namely, even and odd part,
\begin{eqnarray*}
H^2(\Gamma(C)_2)
&=&H^2_{\rm even}(\Gamma(C)_2)\oplus H^2_{\rm odd}(\Gamma(C)_2)\\
&=&\Bigl(\bigoplus_{
                  \lambda\in{\mathcal R}'_{2,\rm even}
                  }
                  \pi_{\lambda}\otimes \pi_{\lambda}^*
   \Bigr)
   \oplus
   \Bigl(\bigoplus_{
                    \lambda\in{\mathcal R}'_{2,\rm odd}}
                  \pi_{\lambda}\otimes \pi_{\lambda}^*
   \Bigr).
\end{eqnarray*}
The even part 
\begin{eqnarray*}
H^2(\Gamma(C)_2)_{\rm even}&=&\bigoplus_{\lambda\in{\mathcal R}'_{2,\rm even}}\pi_{\lambda}\otimes\pi^*_{\lambda}\\
&=&\bigoplus_{\scriptstyle (\lambda_1,\ldots,\lambda_r)\in\Bbb Z^r
              \atop
              \scriptstyle -r>\lambda_1\geq\ldots\geq\lambda_r}\pi_\lambda\otimes\pi_\lambda^*
\end{eqnarray*}
coincides with the Hardy space $H^2(\Gamma(C))$ on the semigroup
$\Gamma(C)$. \\
Our goal in now is to identify the odd part
\begin{eqnarray*}
H^2(\Gamma(C)_2)_{\rm odd}&=&\bigoplus_{\lambda\in{\mathcal R}'_{2,\rm odd}}\pi_{\lambda}\otimes\pi^*_{\lambda}\\
&=&\bigoplus_{\scriptstyle (\lambda_1,\ldots,\lambda_r)\in\Bbb Z^r+{1\over2}
              \atop
              \scriptstyle -(r+{1\over2})\geq\lambda_1\geq\ldots\geq\lambda_r}\pi_\lambda\otimes\pi_\lambda^*
\end{eqnarray*}
 with the classical Hardy space $H^2(Sp(2r, \Bbb R)/U(2r))$.
\begin{theorem}[{\cite{Koufany-Orsted1}}, {\cite{Koufany-Orsted3}}]
The operator $\C_{r+{1\over2}}$ given by (\ref{Cayley}) induces a unitary isomorphism
$$H^2(\Gamma(C)_2)_{\rm odd}\simeq H^2(Sp(2r,\Bbb R)/U(2r)).$$
\end{theorem}

\begin{corollary}[{\cite{Koufany-Orsted1}}, {\cite{Koufany-Orsted3}}]
Under the action of $Mp(r, \Bbb R)\times Mp(r, \Bbb R)$ the Hardy space $H^2(T_\Omega)$ can be decomposed into a direct sum of the `odd' holomorphic discrete series representations of $Mp(r,\Bbb R)$, \ie.
$$H^2(Sp(2r,\Bbb R)/U(2r))_{|_{Mp(r,\Bbb R)\times Mp(r,\Bbb R)}}=\bigoplus_{\lambda\in {\mathcal R}'_{2,\rm odd}}\pi_\lambda\otimes \pi^*_\lambda.$$
\end{corollary}


\begin{corollary}[{\cite{Koufany-Orsted1}}, {\cite{Koufany-Orsted3}}]
Let $K_{\rm odd}$ be the kernel corresponding to $H^2(\Gamma(C)_2)_{\rm odd}$. Then for every $\gamma_1, \gamma_2\in \Gamma(C)_2$
$$K_{\rm odd}(\gamma_1,\gamma_2)=\Det(J-\gamma_2^*J\gamma_1)^{-(r+1/2)}.$$
\end{corollary}


\begin{corollary}[{\cite{Koufany-Orsted1}}, {\cite{Koufany-Orsted3}}]
On the interior of the metaplectic semigroup the distribution
$\Det(I-\gamma)^{-(r+1/2)}$ has the following expansion
$$\Det(I-\gamma)^{-(r+1/2)}=\sum_{\lambda\in{\mathcal R}'_{2,\rm
odd}}d_\lambda\tr(\pi_\lambda(\gamma)),$$
where $d_\lambda$ is the formal dimension of $\pi_\lambda$.
\end{corollary}

\bigskip

The Bergman space \index{Bergman space} on $\Gamma(C)$ is ${\mathcal H}_{2r+1}(\Gamma(C))$ and its
reproducing kernel is given by
$$K_B(\gamma_1,\gamma_2)=\Det(J-\gamma_2^*J\gamma_1)^{-(2r+1)}.$$

\begin{corollary}[{\cite{Koufany-Orsted3}}]
The Bergman kernel $K_B$ on the semigroup $\Gamma(C)$ is the  square of
the odd part $K_{\rm odd}$ of the Cauchy-Szeg{\"o} kernel for $\Gamma(C)_2$.
\end{corollary}

\newpage

\section{The case of $G=SO^*(2l)$}\label{section16}
{\sl We study the case $G=SO^*(2l)$ as we did in  section \ref{section15}. We prove that the classical Hardy space is a proper subspace of odd part
  of the
  Olshanski\u{\i} Hardy space}\\

Let $G=SO^*(2l)$ realized as a subgroup of $U(l,l)$,
$$G=\{g\in {SO}^*(2l,\Bbb C) \tq g^*Jg=J\}, \,\,\,J=\begin{pmatrix}-I_l&0\cr
0&I_l\end{pmatrix}$$
The Hardy parameter in this case is $N/R=l(2l-1)/l=2l-1$ and the
Koecher norm ``$\det$'' is the square root of the usual determinant
``$\Det$'' ($\det=\Det^{1/2}$). Thus the operator $\C_{N\over R}=\C_{2l-1}$,
\begin{equation*}\label{Cayley2}
f=\C_{2l-1}(F)\,\,\,:\,\,\,f(Z)=\Det(Z+iJ)^{-(l-1/2)}F(\gamma)
\end{equation*}
provides an equivariant embedding of the classical Hardy space\\ 
$H^2(SO^*(4l)/U(2l))$ into the
odd part of the Hardy space $H^2(\Gamma(C)_2)_{\rm odd}$ on the double
covering semigroup $\Gamma(C)_2$ of the minimal semigroup
$$\Gamma(C)=\{\gamma\in{SO}^*(2l,\Bbb C) \tq J-\gamma^*J\gamma\geq0\}$$
(because of the square root in $\Det(Z+iJ)^{(l-1/2)}$).
This is just like the symplectic case. We will identify the maximal compact subgroup with
$U(l)$ as in the above section. The determinant factor is again exactly the Jacobian to
 a power such that we have preservation of ${\rm L}^2-$norms on the respective boundaries. Then $\frak t_{\Bbb
R}$ is given by
the same formula as in $Sp(r,\Bbb R)$ case. Let $\epsilon_1, \epsilon_2, \ldots, \epsilon_l$ be the
canonical basis of $\frak t^*_{\Bbb R}=\frak t_{\Bbb R}=\Bbb R^l$. The
root system $\Delta=\Delta(\frak g_{\Bbb C},\frak t_{\Bbb C})$ is of type $D_l$ :
\begin{eqnarray*}
\Delta&=&\{\pm\epsilon_i\pm\epsilon_j \tq 1\leq i<j\leq l\},\\
\Delta^+&=&\{\epsilon_i\pm\epsilon_j \tq 1\leq i<j\leq l\},\\
\Delta^+_{\frak k}&=&\{\epsilon_i-\epsilon_j \tq 1\leq i<j\leq l\},\\
\Delta^+_{\frak p}&=&\{\epsilon_i+\epsilon_j \tq 1\leq i<j\leq l\},\\
\rho&=&(l-1)\epsilon_1+(l-2)\epsilon_2+\ldots+\epsilon_{l-1}\\
&\simeq&(l-1,l-2,\ldots,1,0).
\end{eqnarray*}
The set of highest weights relative to the positive roots of $SO^*(2l)$ is 
$${\mathcal R}=\{\lambda=(\lambda_1\,\ldots,\lambda_l)\in \Bbb Z^l \tq
\lambda_1\geq\ldots\geq\lambda_l\},$$
and $\lambda\in{\mathcal R}$ satisfies to the Harish-Chandra condition if and only if
\begin{equation*}
-2l+3>\lambda_1+\lambda_2.
\end{equation*}
Therefore, the odd holomorphic discrete series representations of the
double covering group $G_2$ of ${SO}^*(2l)$ are those irreducible
unitary representations $\pi_\lambda$, square-integrable with a
highest weight $\lambda=(\lambda_1,\ldots,\lambda_l)\in\Bbb
Z^l-{1\over2}$ such that $0\geq\lambda_1\geq\ldots\geq\lambda_l$ and satisfying
\begin{equation*}
-2l+2\geq\lambda_1+\lambda_2.
\end{equation*}
Let  ${\mathcal R}'_{2,\rm odd}$ denotes the set of these $\lambda$'s.
The Hardy space $H^2(\Gamma(C)_2)$ on the minimal cone $\Gamma(C)_2$
has then the following decomposition
\begin{equation*}
H^2(\Gamma(C)_2)_{\rm odd}=\bigoplus_{\lambda\in{\mathcal R}'_{2,\rm
odd}}\pi_\lambda\otimes{\pi}^*_\lambda.
\end{equation*}

\begin{theorem}[{\cite{Koufany-Orsted3}}]
The classical Hardy space $H^2(SO^*(4l)/U(2l))$ is a proper invariant
subspace of the ``non--classical'' Hardy space $H^2(\Gamma(C)_2)_{\rm
odd}$.
\end{theorem}

\begin{corollary}[{\cite{Koufany-Orsted3}}]
The representation of $SO^*(2l)\times SO^*(2l)$ in the Hardy space
$H^2(\Gamma(C))$ cannot be obtained by a restriction of a
representation of the holomorphic discrete series of $SO^*(4l)$ nor
any continuation of this, such as the Hardy space.
\end{corollary}

Let ${\mathcal H}^2(\Gamma(C))$ the conformal image of $H^2(SO^*(4l)/U(2l))$
via the operator $\C_{2l-1}$.

\begin{corollary}[{\cite{Koufany-Orsted3}}]
${\mathcal H}^2(\Gamma(C))$ is a reproducing kernel Hilbert space and its
reproducing kernel $K$ is the pre-image of the Cauchy--Szeg{\"o}
kernel of $H^2(SO^*(4l)/U(2l))$, \ie.
$$K(\gamma_1,\gamma_2)=\Det(J-\gamma_2^*J\gamma_1)^{-(l-1/2)}.$$
\end{corollary}

\begin{corollary}[{\cite{Koufany-Orsted3}}]
On $\Gamma(C)^\circ$  the holomorphic function $\Det(I-\gamma)^{-(l-1/2)}$ has the
following expansion
$$\Det(I-\gamma)^{-(l-1/2)}=\sum_{-l+1/2\geq\lambda_1\geq\ldots\geq\lambda_l}d_\lambda\tr(\pi_\lambda(\gamma)),$$
where $d_\lambda$ is the formal dimension of $\pi_\lambda$. 
\end{corollary}

\newpage

\section{The case of $G=U(p,q)$}\label{section17}
{\it We study the case $G=U(p,q)$ as we did in  section
  \ref{section15}. Here we do not need double covering semigroup.
  We prove that the classical Hardy space is a proper subspace of
  of the
  Olshanski\u{\i} Hardy space}\\

Wee fix $G=U(p,q)$  realized by 
$$G=U(p,q)=\{g\in GL(n,\Bbb C) \tq
g^*Jg=J\},\,J=\begin{pmatrix}-I_p&0\cr0&I_q\end{pmatrix}$$
where $n=p+q$. In this case the Hardy parameter is ${N/R}=n^2/n=n$ and the Koecher
norm ``$\det$'' is the usual determinant ``Det''. Therefore the operator
$\C_{N\over R}$ given by
\begin{equation*}
f=\C_n(F) \,\, : \,\, f(Z)=\Det(Z+iJ)^{-n}F(\gamma)
\end{equation*}
may be an intertwining operator between the classical Hardy space
$H^2(SU(n,n)/S\bigl(U(n)\times U(n)\bigr)$ and the Hardy space
$H^2(\Gamma(C))$ over the semigroup
$$\Gamma(C)=\{\gamma\in GL(n,\Bbb C)\tq J-\gamma^*J\gamma\geq0 \}.$$
To study unitary representations of $G=U(p,q)$ we identify it with $(U(1)\times SU(p,q))/{\Bbb
Z_{p+q}}$. Thus the unitary irreducible representations of $G$ are those
of $U(1)\times SU(p,q)$ that are trivial on $(\zeta,\zeta^{-1}I_n)$ as in Remark \ref{remark2.1}, where $I_n$ is
the identity matrix and $\zeta^n = 1$. Therefore, the holomorphic discrete
series representations of $G$ that we are interested in are 
\begin{equation*}
\pi_{\lambda,k}({e}^{i\theta}g)={e}^{i k\theta}\pi_{\lambda}(g) \,\, , \,\,
g\in SU(p,q),\,\,\theta\in\Bbb R
\end{equation*}
where $k\in\Bbb Z$ and $\pi_\lambda$ are the holomorphic discrete
series representations of $SU(p,q)$, realized on ${\mathcal
D}=SU(p,q)/S\bigl(U(p)\times U(q)\bigr)$, for example in the scalar case: 
$$\bigl(\pi_\lambda(g)f\bigr)(Z)=\Det(CZ+D)^{-\lambda}f\bigl((AZ+B)(CZ+D)^{-1}\bigr), \,\,\, g^{-1}=\begin{pmatrix}A&B\cr C&D\end{pmatrix},$$
with $\lambda$ an integer, and in general $\lambda=(\lambda_1, \ldots,
\lambda_n)\in \Bbb Z^n$.
There will an underlying parity condition to make the representations trivial on
$\Bbb Z_n$ as above; for example in the scalar case we must have 
that $k-q\lambda$ is divisible by $n$.\\
Let $\frak t\subset \frak k$ be a Cartan subalgebra consisting of
diagonal matrices with purely imaginary values and $\frak t_{\Bbb
R}=i\frak t$. Then the root system $\Delta=\Delta(\frak g_{\Bbb
C},\frak t_{\Bbb C})$ is of type $A_{n-1}$ :
\begin{eqnarray*}
\Delta&=&\{\epsilon_i-\epsilon_j\tq  1\leq i\not=j \leq n\},\\
\Delta^+&=&\{\epsilon_i-\epsilon_j \tq 1\leq i<j\leq n \},\\
\Delta^+_{\frak k}&=&\{\epsilon_i-\epsilon_j \tq 1\leq i<j\leq p
\,\,\,{\rm or}\,\,\,
p+1\leq i<j\leq n \},\\
\Delta^+_{\frak p}&=&\{\epsilon_i-\epsilon_j \tq 1\leq i\leq p \,\,\,{\rm and}\,\,\,
p+1\leq j\leq n \},\\
2\rho&=&(n-1)\epsilon_1+(n-3)\epsilon_2+\ldots-(n-3)\epsilon_{n-1}-(n-1)\epsilon_n\cr
&\simeq&(n-1, n-3, n-5, \ldots, -n+3, -n+1),
\end{eqnarray*}
where $\epsilon_1, \ldots, \epsilon_n$ is the canonical basis of
$\frak t_{\Bbb R}^*\simeq\frak t_{\Bbb R}\simeq\Bbb R^n$. Then the holomorphic discrete series representations of $SU(p,q)$ are the above representations $\pi_\lambda$  with $\lambda=(\lambda_1,\ldots,\lambda_n)\in \Bbb Z^n$ satisfying
$$\lambda_i-\lambda_{i+1}\geq 0, \;\; i\not=p, \;\; 1\leq i\leq n-1,$$ 
and the
Harish-Chandra condition
$$\lambda_n-\lambda_1>n-1.$$
Let $\pi_{\lambda,k}$ be an irreducible unitary
representation of $G$ with highest weight
$(\lambda,k)$, $\lambda=(\lambda_1, \ldots, \lambda_n)$. 
 Then  we prove that $\pi_{\lambda,k}$ is a $C-$dissipative representation of the holomorphic discrete series if and only if $(\lambda,k)$ belongs to the set ${\mathcal R}_{\rm diss}$ of $(\lambda_1, \ldots ,\lambda_n,k)\in\Bbb
Z^{n+1}$ such that
$$
\left\{\begin{matrix}
\lambda_n-\lambda_1>n-1\\
0\geq\lambda_1\geq\ldots\geq\lambda_p, \; \lambda_{p+1}\geq\ldots\geq\lambda_n\geq0\\
[\lambda]-n\lambda_n\leq k\leq [\lambda]-n\lambda_1.\\
\end{matrix}\right.
$$
Hence the Hardy space on the semigroup $\Gamma(C)$ has the following
decomposition 
\begin{equation*}
H^2(\Gamma(C))=\bigoplus_{(\lambda,k)\in{\mathcal R}_{\rm
diss}}\pi_{\lambda,k}\otimes\pi^*_{\lambda,k}.
\end{equation*}

\begin{theorem}[{\cite{Koufany-Orsted3}}]
The classical Hardy space $H^2\bigl(SU(n,n)/S\bigl(U(n)\times
U(n)\bigr)\bigr)$ is a proper invariant subspace of the ``non--classical'' Hardy space $H^2(\Gamma(C))$.
\end{theorem}

\begin{corollary}[{\cite{Koufany-Orsted3}}]
The representation of $S\bigl(U(p,q)\times U(p,q)\bigr)$ in the Hardy
space $H^2(\Gamma(C))$ cannot be obtained by a restriction of a
representation of the holomorphic discrete series of $SU(n,n)$ nor any
analytic continuation of this, such as the Hardy space.
\end{corollary}

Let ${\mathcal H}^2(\Gamma(C))$ be the conformal image of $H^2\bigl(SU(n,n)/S\bigl(U(n)\times
U(n)\bigr)\bigr)$ via the operator $\C_n$.

\begin{corollary}[{\cite{Koufany-Orsted3}}]
${\mathcal H}^2(\Gamma(C))$
is a reproducing Hilbert space and its reproducing kernel $K$ is the
pre-image of the Cauchy-Szeg{\"o} kernel of $H^2\bigl(SU(n,n)/S\bigl(U(n)\times
U(n)\bigr)\bigr)$, \ie.
$$K(\gamma_1,\gamma_2)=\Det(J-\gamma_2^*J\gamma_1)^{-n}.$$
\end{corollary}

\begin{corollary}[{\cite{Koufany-Orsted3}}]
On $\Gamma(C)^\circ$ the holomorphic function $\Det(I-\gamma)^{-n}$ has the following expansion
$$\Det(I-\gamma)^{-n}=\sum_{(\lambda,k)\in {\mathcal R}_{\rm decay}}
d_{\lambda,k}\tr(\pi_{\lambda,k}(\gamma)),$$
where $d_{\lambda,k}$ is the formal dimension of $\pi_{\lambda,k}$
and 
${\mathcal R}_{\rm decay}$ is the set of $(\lambda_1, \ldots ,\lambda_n,k)\in\Bbb
Z^{n+1}$ such that
$$
\left\{\begin{matrix}
\lambda_n-\lambda_1>n-1\cr
0\geq\lambda_1\geq\ldots\geq\lambda_p, \; \lambda_{p+1}\geq\ldots\geq\lambda_n\geq0\cr
[\lambda]-n(\lambda_n+n)\leq k\leq[\lambda]-n(\lambda_1+n)\cr
\end{matrix}\right.
$$
\end{corollary}

\printindex

\end{document}